\documentstyle[amssymb]{amsart}
\def\grayrule{\special{ps:: currentgray 0.9 setgray}
\hbox to 0pt{\hss
\vrule height 0.3 in depth 0.2 in width 0.49 in \hss}
\special{ps:: setgray}}

\newcommand{\ZFCa}{\text{\normalshape\sf ZFC}}
\newcommand{\CH}{\text{\normalshape\sf CH}}
\newcommand{\coll}{\text{\normalshape\sf Coll}}

\newcommand{\rest}{{\mathord{\restriction}}}

\newcommand{\dom}{{\text{\normalshape\sf {dom}}}}

\newcommand{\spli}{{\text{\normalshape\sf{split}}}}

\newcommand{\supp}{{\text{\normalshape\sf {supp}}}}
\newcommand{\suc}{{\text{\normalshape\sf {succ}}}}

\newcommand{\Seq}{{\text{\normalshape\sf {Seq}}}}
\newcommand{\lh}{{\text{\normalshape\sf {lh}}}}
\newcommand{\Con}{{\text{\normalshape\sf {Con}}}}
\newcommand{\QED}{\hspace{0.1in} \square \vspace{0.1in}}
\newcommand{\forces}{\mathrel{\|}\joinrel\mathrel{-}}
\newtheorem{theorem}{Theorem}[section]
\newtheorem{lemma}[theorem]{Lemma}
\newtheorem{claim}[theorem]{Claim}

\newcommand{\Proof}{{\sc Proof} \hspace{0.2in}}
\newcommand{\lft}[2]{\mathopen\ifcase#1{}\oo\or
                        \big#2\or\Big#2\else\oo\fi} 
\newcommand{\rgt}[2]{\mathclose\ifcase#1{}\oo\or
                        \big#2\or\Big#2\else\oo\fi} 
\theoremstyle{definition}
\newtheorem{definition}[theorem]{Definition}
\begin{document}

\title{The Cicho\'{n} Diagram}
\author{Tomek Bartoszy\'{n}ski}
\address{Department of Mathematics 
University of California
Berkeley, CA  94720}
\curraddr{Department of Mathematics  
Boise State University, Boise, Idaho 83725}
\email{tomek@@math.idbsu.edu}
\author{Haim  Judah} \thanks{The second author
would like to thank NSF under grant DMS-8505550 and MSRI for partial
support.} 
\address{ Department of Mathematics 
University of California, Berkeley 
and 
Mathematical Sciences Research Institute, 1000 Centennial Drive,
Berkeley, California}
\curraddr{Department of Mathematics 
Bar Ilan University 
Ramat Gan, Israel 52900}
\email{judah@@bimacs.cs.biu.ac.il}
\author{Saharon Shelah} 
\thanks{The third author would like to thank 
US-Israel Binational Science Foundation and MSRI for partial support.
Publication number 368}
\address{Department of Mathematics 
Rutgers University, New Brunswick, New Jersey and Mathematical
Sciences Research Institute, 1000 Centennial Drive, 
Berkeley, California and Department
of Mathematics  Hebrew University Jerusalem}

\curraddr{Department
of Mathematics,  Hebrew University, Jerusalem}
\email{shelah@@sunrise.huji.ac.il}

\maketitle
\begin{abstract} We conclude the discussion of
additivity, Baire number, uniformity and covering for measure and category
by constructing the remaining
5 models. Thus we complete the analysis of Cicho\'{n}'s diagram.
\end{abstract}
\section{Introduction}
The goal of this paper is to  describe the relationship between basic
properties  
of
measure and category.

\begin{definition}
Let ${\cal N}$ and ${\cal M}$ denote the
ideals of null subsets of the real line and meager
subsets of the real line respectively.

Define the following ten sentences:

${\bold A}(m)$ $\equiv$ unions of fewer than $2^{\boldsymbol\aleph_{0}}$ null sets is
null, 

${\bold B}(m)$ $\equiv$ $\Re$ is not the union of fewer than
$2^{\boldsymbol\aleph_{0}}$ 
null sets,

${\bold U}(m)$ $\equiv$ every subset of $\Re$ of size less than
$2^{\boldsymbol\aleph_{0}}$ 
is null,

${\bold C}(m)$ $\equiv$ ideal of null sets does not have a basis
of size less than
$2^{\boldsymbol\aleph_{0}}$ .

Sentences ${\bold A}(c),\ {\bold B}(c),\ {\bold U}(c)$ and ${\bold C}(c)$ are 
defined
analogously by replacing word ``null'' by the word ``meager'' in the
definitions 
above.

In addition define

$w{\bold D} \ \equiv \ \forall F \subset [\omega^{\omega}]^{<2^{\boldsymbol\aleph_{0}}}\
\exists g \in \omega^{\omega} \ \forall f \in F \ \exists^{\infty}n \
f(n)<g(n)$

and

${\bold D} \ \equiv \ \forall F \subset [\omega^{\omega}]^{<2^{\boldsymbol\aleph_{0}}}\
\exists g \in \omega^{\omega} \ \forall f \in F \ \forall^{\infty}n \
f(n)<g(n)$.
\end{definition}

The relationship between these sentences is
described in the following diagram which is called Cicho\'{n}'s diagram:
$$\begin{array}{ccccccc}
{\bold B}(m) &
\rightarrow &
{\bold U}(c) &
\rightarrow &
{\bold C}(c)&
\rightarrow &
{\bold C}(m) \\
\ &
\ &
\uparrow &
\ &
\uparrow &
\ &
\  \\
\smash{\bigg\uparrow} &
\ &
{\bold D} &
\rightarrow  &
w{\bold D} &
\ &
\smash{\bigg\uparrow} \\
\ &
\ &
\uparrow &
\ &
\uparrow &
\ &
\ \\
{\bold A}(m) &
\rightarrow &
{\bold A}(c) &
\rightarrow &
{\bold B}(c) &
\rightarrow &
{\bold U}(m)
\end{array}$$

In addition
$${\bold A}(c) \equiv {\bold B}(c) \ \& \ {\bold D}$$
and
$${\bold C}(c) \equiv {\bold U}(c) \ \vee \ w{\bold D}  .$$

The proofs of these inequalities can be found in \cite{Ba}, \cite{F} and \cite{M1}.

In context of this diagram a natural question arises: 

  Are those the only implications
between these sentences that are provable in
\ZFCa ?

It turns out that the answer to this question is positive. Every
combination 
of those sentences which does not contradict the implications in the
diagram 
is consistent with \ZFCa.
This is proved in step-by-step fashion and this paper contains
constructions 
of the last 5 models. 

The tables below contain all known results on the subject.
They are not symmetric but still one can recognize
some patterns here. 
Let ${\cal L}$ be the set of sentences obtained
from sentences
${\bold A}$,
${\bold B}$,
${\bold U}$,
${\bold C}$, ${\bold D}$ and $w{\bold D}$ using logical connectives.
Define $^{\star} : {\cal L} \longrightarrow {\cal L}$ as\\

$\varphi^{\star} = \left\{ \begin{array}{ll}
\neg \psi^{\star} & \hbox{if } \varphi=\neg \psi \\
\psi_{1}^{\star} \vee \psi_{2}^{\star} & \hbox{if } \varphi=\psi_{1} \vee \psi_{2}\\
\neg {\bold C} & \hbox{if } \varphi = {\bold A} \\
\neg {\bold U} & \hbox{if } \varphi = {\bold B} \\
\neg {\bold B} & \hbox{if } \varphi = {\bold U} \\
\neg {\bold A} & \hbox{if } \varphi = {\bold C} \\
\neg w{\bold D} & \hbox{if } \varphi = {\bold D} \\
\neg {\bold D} & \hbox{if } \varphi = w{\bold D}
\end{array} \right.$
for $\varphi \in {\cal L}$. 

It turns out that if $\varphi$ is consistent with $\ZFCa$ then
$\varphi^{\star} $ is consistent with $\ZFCa$.
Moreover, in most cases one can find a notion of forcing ${\bold P}$ such
that $\omega_{2}$-iteration of ${\bold P}$ over a model for $\CH$ gives
a model for $\varphi$ while
$\omega_{1}$-iteration of ${\bold P}$ over a model for ${\bold {MA}} \ \&  \ \neg \CH$
gives a model for $\varphi^{\star}$.

The first table known as, the Kunen-Miller chart, gives consistency
results 
concerning sentences 
${\bold A}$,
${\bold B}$,
${\bold U}$,
${\bold C}$ only. It was completed by H. Judah and S. Shelah in \cite{JS}.
The remaining three tables give corresponding information including
all 3 consistent combinations of $\bold D$ and $w{\bold D}$.

\begin{tabular}{|c|c|c|c|c|c|c|c|c|c|}\hline
\multicolumn{3}{|l|}
{\makebox[0.35in]{\rule[-0.2in]{0in}{0.5in}
\rule[-0.2in]{0in}{0.5in}}}&
\makebox[0.35in]{\rule[-0.2in]{0in}{0.5in}
{\small Add}
\rule[-0.2in]{0in}{0.5in}}&
\makebox[0.35in]{\rule[-0.2in]{0in}{0.5in}
${\bold T}$
\rule[-0.2in]{0in}{0.5in}}&
\makebox[0.35in]{\rule[-0.2in]{0in}{0.5in}
${\bold F}$
\rule[-0.2in]{0in}{0.5in}}&
\makebox[0.35in]{\rule[-0.2in]{0in}{0.5in}
${\bold F}$
\rule[-0.2in]{0in}{0.5in}}&
\makebox[0.35in]{\rule[-0.2in]{0in}{0.5in}
${\bold F}$
\rule[-0.2in]{0in}{0.5in}}&
\makebox[0.35in]{\rule[-0.2in]{0in}{0.5in}
${\bold F}$
\rule[-0.2in]{0in}{0.5in}}&
\makebox[0.35in]{\rule[-0.2in]{0in}{0.5in}
${\bold F}$
\rule[-0.2in]{0in}{0.5in}}\\
\cline{4-10}
\multicolumn{3}{|r|}
{\makebox[0.35in]{\rule[-0.2in]{0in}{0.5in}
{\small Category}
\rule[-0.2in]{0in}{0.5in}}}&
\makebox[0.35in]{\rule[-0.2in]{0in}{0.5in}
{\small Baire}
\rule[-0.2in]{0in}{0.5in}}&
\makebox[0.35in]{\rule[-0.2in]{0in}{0.5in}
${\bold T}$
\rule[-0.2in]{0in}{0.5in}}&
\makebox[0.35in]{\rule[-0.2in]{0in}{0.5in}
${\bold T}$
\rule[-0.2in]{0in}{0.5in}}&
\makebox[0.35in]{\rule[-0.2in]{0in}{0.5in}
${\bold T}$
\rule[-0.2in]{0in}{0.5in}}&
\makebox[0.35in]{\rule[-0.2in]{0in}{0.5in}
${\bold F}$
\rule[-0.2in]{0in}{0.5in}}&
\makebox[0.35in]{\rule[-0.2in]{0in}{0.5in}
${\bold F}$
\rule[-0.2in]{0in}{0.5in}}&
\makebox[0.35in]{\rule[-0.2in]{0in}{0.5in}
${\bold F}$
\rule[-0.2in]{0in}{0.5in}}\\ \cline{4-10}
\multicolumn{3}{|c|}
{\makebox[0.35in]{\rule[-0.2in]{0in}{0.5in}
{\small Measure}
\rule[-0.2in]{0in}{0.5in}}}&
\makebox[0.35in]{\rule[-0.2in]{0in}{0.5in}
{\small Unif}
\rule[-0.2in]{0in}{0.5in}}&
\makebox[0.35in]{\rule[-0.2in]{0in}{0.5in}
${\bold T}$
\rule[-0.2in]{0in}{0.5in}}&
\makebox[0.35in]{\rule[-0.2in]{0in}{0.5in}
${\bold T}$
\rule[-0.2in]{0in}{0.5in}}&
\makebox[0.35in]{\rule[-0.2in]{0in}{0.5in}
${\bold F}$
\rule[-0.2in]{0in}{0.5in}}&
\makebox[0.35in]{\rule[-0.2in]{0in}{0.5in}
${\bold T}$
\rule[-0.2in]{0in}{0.5in}}&
\makebox[0.35in]{\rule[-0.2in]{0in}{0.5in}
${\bold F}$
\rule[-0.2in]{0in}{0.5in}}&
\makebox[0.35in]{\rule[-0.2in]{0in}{0.5in}
${\bold F}$
\rule[-0.2in]{0in}{0.5in}}\\ \hline
\makebox[0.35in]{\rule[-0.2in]{0in}{0.5in}
{\small Add}
\rule[-0.2in]{0in}{0.5in}}&
\makebox[0.35in]{\rule[-0.2in]{0in}{0.5in}
{\small Baire}
\rule[-0.2in]{0in}{0.5in}}&
\makebox[0.35in]{\rule[-0.2in]{0in}{0.5in}
{\small Unif}
\rule[-0.2in]{0in}{0.5in}}&
\makebox[0.35in]{\rule[-0.2in]{0in}{0.5in}
{\small Cov}
\rule[-0.2in]{0in}{0.5in}}&
\makebox[0.35in]{\rule[-0.2in]{0in}{0.5in}
${\bold T}$
\rule[-0.2in]{0in}{0.5in}}&
\makebox[0.35in]{\rule[-0.2in]{0in}{0.5in}
${\bold T}$
\rule[-0.2in]{0in}{0.5in}}&
\makebox[0.35in]{\rule[-0.2in]{0in}{0.5in}
${\bold T}$
\rule[-0.2in]{0in}{0.5in}}&
\makebox[0.35in]{\rule[-0.2in]{0in}{0.5in}
${\bold T}$
\rule[-0.2in]{0in}{0.5in}}&
\makebox[0.35in]{\rule[-0.2in]{0in}{0.5in}
${\bold T}$
\rule[-0.2in]{0in}{0.5in}}&
\makebox[0.35in]{\rule[-0.2in]{0in}{0.5in}
${\bold F}$
\rule[-0.2in]{0in}{0.5in}}\\ \hline
\makebox[0.35in]{\rule[-0.2in]{0in}{0.5in}
${\bold T}$
\rule[-0.2in]{0in}{0.5in}}&
\makebox[0.35in]{\rule[-0.2in]{0in}{0.5in}
${\bold T}$
\rule[-0.2in]{0in}{0.5in}}&
\makebox[0.35in]{\rule[-0.2in]{0in}{0.5in}
${\bold T}$
\rule[-0.2in]{0in}{0.5in}}&
\makebox[0.35in]{\rule[-0.2in]{0in}{0.5in}
${\bold T}$
\rule[-0.2in]{0in}{0.5in}}& \makebox[0.35in]{\rule[-0.2in]{0in}{0.5in}
$A$
\rule[-0.2in]{0in}{0.5in}}&
\makebox[0.35in]{\rule[-0.2in]{0in}{0.5in}
\grayrule
\rule[-0.2in]{0in}{0.5in}}&
\makebox[0.35in]{\rule[-0.2in]{0in}{0.5in}
\grayrule
\rule[-0.2in]{0in}{0.5in}}&
\makebox[0.35in]{\rule[-0.2in]{0in}{0.5in}
\grayrule
\rule[-0.2in]{0in}{0.5in}}&
\makebox[0.35in]{\rule[-0.2in]{0in}{0.5in}
\grayrule
\rule[-0.2in]{0in}{0.5in}}&
\makebox[0.35in]{\rule[-0.2in]{0in}{0.5in}
\grayrule
\rule[-0.2in]{0in}{0.5in}}\\ \hline
\makebox[0.35in]{\rule[-0.2in]{0in}{0.5in}
${\bold F}$
\rule[-0.2in]{0in}{0.5in}}&
\makebox[0.35in]{\rule[-0.2in]{0in}{0.5in}
${\bold T}$
\rule[-0.2in]{0in}{0.5in}}&
\makebox[0.35in]{\rule[-0.2in]{0in}{0.5in}
${\bold T}$
\rule[-0.2in]{0in}{0.5in}}&
\makebox[0.35in]{\rule[-0.2in]{0in}{0.5in}
${\bold T}$
\rule[-0.2in]{0in}{0.5in}}&
\makebox[0.35in]{\rule[-0.2in]{0in}{0.5in}
$  B$
\rule[-0.2in]{0in}{0.5in}}&
\makebox[0.35in]{\rule[-0.2in]{0in}{0.5in}
$  C$
\rule[-0.2in]{0in}{0.5in}}&
\makebox[0.35in]{\rule[-0.2in]{0in}{0.5in}
\grayrule
\rule[-0.2in]{0in}{0.5in}}&
\makebox[0.35in]{\rule[-0.2in]{0in}{0.5in}
$  D$
\rule[-0.2in]{0in}{0.5in}}&
\makebox[0.35in]{\rule[-0.2in]{0in}{0.5in}
\grayrule
\rule[-0.2in]{0in}{0.5in}}&
\makebox[0.35in]{\rule[-0.2in]{0in}{0.5in}
\grayrule
\rule[-0.2in]{0in}{0.5in}}\\ \hline
\makebox[0.35in]{\rule[-0.2in]{0in}{0.5in}
${\bold F}$
\rule[-0.2in]{0in}{0.5in}}&
\makebox[0.35in]{\rule[-0.2in]{0in}{0.5in}
${\bold T}$
\rule[-0.2in]{0in}{0.5in}}&
\makebox[0.35in]{\rule[-0.2in]{0in}{0.5in}
${\bold F}$
\rule[-0.2in]{0in}{0.5in}}&
\makebox[0.35in]{\rule[-0.2in]{0in}{0.5in}
${\bold T}$
\rule[-0.2in]{0in}{0.5in}}&
\makebox[0.35in]{\rule[-0.2in]{0in}{0.5in}
\grayrule
\rule[-0.2in]{0in}{0.5in}}&
\makebox[0.35in]{\rule[-0.2in]{0in}{0.5in}
\grayrule
\rule[-0.2in]{0in}{0.5in}}&
\makebox[0.35in]{\rule[-0.2in]{0in}{0.5in}
\grayrule
\rule[-0.2in]{0in}{0.5in}}&
\makebox[0.35in]{\rule[-0.2in]{0in}{0.5in}
$  E\!=\!E^\star$
\rule[-0.2in]{0in}{0.5in}}&
\makebox[0.35in]{\rule[-0.2in]{0in}{0.5in}
\grayrule
\rule[-0.2in]{0in}{0.5in}}&
\makebox[0.35in]{\rule[-0.2in]{0in}{0.5in}
\grayrule
\rule[-0.2in]{0in}{0.5in}}\\ \hline
\makebox[0.35in]{\rule[-0.2in]{0in}{0.5in}
${\bold F}$
\rule[-0.2in]{0in}{0.5in}}&
\makebox[0.35in]{\rule[-0.2in]{0in}{0.5in}
${\bold F}$
\rule[-0.2in]{0in}{0.5in}}&
\makebox[0.35in]{\rule[-0.2in]{0in}{0.5in}
${\bold T}$
\rule[-0.2in]{0in}{0.5in}}&
\makebox[0.35in]{\rule[-0.2in]{0in}{0.5in}
${\bold T}$
\rule[-0.2in]{0in}{0.5in}}&
\makebox[0.35in]{\rule[-0.2in]{0in}{0.5in}
$  F$
\rule[-0.2in]{0in}{0.5in}}&
\makebox[0.35in]{\rule[-0.2in]{0in}{0.5in}
$  G$
\rule[-0.2in]{0in}{0.5in}}&
\makebox[0.35in]{\rule[-0.2in]{0in}{0.5in}
$  H\!=\!H^\star$
\rule[-0.2in]{0in}{0.5in}}&
\makebox[0.35in]{\rule[-0.2in]{0in}{0.5in}
$  I\!=\!I^\star$
\rule[-0.2in]{0in}{0.5in}}&
\makebox[0.35in]{\rule[-0.2in]{0in}{0.5in}
$  G^\star$
\rule[-0.2in]{0in}{0.5in}}&
\makebox[0.35in]{\rule[-0.2in]{0in}{0.5in}
$  F^\star$
\rule[-0.2in]{0in}{0.5in}}\\ \hline
\makebox[0.35in]{\rule[-0.2in]{0in}{0.5in}
${\bold F}$
\rule[-0.2in]{0in}{0.5in}}&
\makebox[0.35in]{\rule[-0.2in]{0in}{0.5in}
${\bold F}$
\rule[-0.2in]{0in}{0.5in}}&
\makebox[0.35in]{\rule[-0.2in]{0in}{0.5in}
${\bold F}$
\rule[-0.2in]{0in}{0.5in}}&
\makebox[0.35in]{\rule[-0.2in]{0in}{0.5in}
${\bold T}$
\rule[-0.2in]{0in}{0.5in}}&
\makebox[0.35in]{\rule[-0.2in]{0in}{0.5in}
\grayrule
\rule[-0.2in]{0in}{0.5in}}&
\makebox[0.35in]{\rule[-0.2in]{0in}{0.5in}
\grayrule
\rule[-0.2in]{0in}{0.5in}}&
\makebox[0.35in]{\rule[-0.2in]{0in}{0.5in}
\grayrule
\rule[-0.2in]{0in}{0.5in}}&
\makebox[0.35in]{\rule[-0.2in]{0in}{0.5in}
$  D^\star$
\rule[-0.2in]{0in}{0.5in}}&
\makebox[0.35in]{\rule[-0.2in]{0in}{0.5in}
$  C^\star$
\rule[-0.2in]{0in}{0.5in}}&
\makebox[0.35in]{\rule[-0.2in]{0in}{0.5in}
$  B^\star$
\rule[-0.2in]{0in}{0.5in}}\\ \hline
\makebox[0.35in]{\rule[-0.2in]{0in}{0.5in}
${\bold F}$
\rule[-0.2in]{0in}{0.5in}}&
\makebox[0.35in]{\rule[-0.2in]{0in}{0.5in}
${\bold F}$
\rule[-0.2in]{0in}{0.5in}}&
\makebox[0.35in]{\rule[-0.2in]{0in}{0.5in}
${\bold F}$
\rule[-0.2in]{0in}{0.5in}}&
\makebox[0.35in]{\rule[-0.2in]{0in}{0.5in}
${\bold F}$
\rule[-0.2in]{0in}{0.5in}}&
\makebox[0.35in]{\rule[-0.2in]{0in}{0.5in}
\grayrule
\rule[-0.2in]{0in}{0.5in}}&
\makebox[0.35in]{\rule[-0.2in]{0in}{0.5in}
\grayrule
\rule[-0.2in]{0in}{0.5in}}&
\makebox[0.35in]{\rule[-0.2in]{0in}{0.5in}
\grayrule
\rule[-0.2in]{0in}{0.5in}}&
\makebox[0.35in]{\rule[-0.2in]{0in}{0.5in}
\grayrule
\rule[-0.2in]{0in}{0.5in}}&
\makebox[0.35in]{\rule[-0.2in]{0in}{0.5in}
\grayrule
\rule[-0.2in]{0in}{0.5in}}&
\makebox[0.35in]{\rule[-0.2in]{0in}{0.5in}
$A^\star$
\rule[-0.2in]{0in}{0.5in}}\\ \hline
\end{tabular} 

\begin{enumerate}

\item[$A$]
$\omega_{2}$-iteration with finite (countable) support of
amoeba reals over a model for $\CH$ or
any model for $\CH$ or
${\bold {MA}}$ works.

\item[$A^\star$]
$\omega_2$-iteration with finite (countable) support of
amoeba reals over a model for $\neg \CH$ or
$\omega_{2}$-iteration of Sacks or Silver reals over
a model for $\CH$.
\item[$ B$] $\omega_{2}$-iteration of random and dominating reals over
  a model 
for $\CH$. \cite{M1}
\item[$ B^\star$] $\omega_{1}$-iteration of random and dominating
  reals over a model
for $\neg \CH \ \& \ {\bold B}(c)$. 
\item[$ C$] $\omega_{2}$-iteration with finite support of random reals
  over a model for $\CH$. \cite{M1}
\item[$ C^\star$] $\omega_{1}$-iteration with finite support of random
  reals  over a model for $\neg \CH \& \ {\bold D}$. \cite{M1}
\item[$  D$] Countable support $\omega_{2}$-iteration of infinitely
  equal reals (see section 3) and   random reals over a model for
  $\CH$. \cite{M1} 
\item[$ D^\star$] $\omega_{2}$-iteration of Laver reals (\cite{JS}). We do not
know if there exists a  notion of forcing ${\bold P}$ such that
$\omega_{2}$-iteration of ${\bold P}$ over a model for $\CH$ gives
$ D$ and 
$\omega_{1}$-iteration of ${\bold P}$ over a model for ${\bold {MA}} \
\& \ \neg \CH$  gives $ D^\star$. 
\item[$ E\!=\!E^\star$] $\boldsymbol\aleph_{2}$ random reals over a
  model for $\CH$. This model is self-dual. 
\item[$ F$] $\omega_{2}$-iteration with finite support of
any $\sigma$-centered notion of forcing adding dominating reals
over a model
for $\CH$ . \cite{M1}
\item[$ F^\star$] $\omega_{1}$-iteration with finite support of
any $\sigma$-centered notion of forcing adding dominating reals
over a model
for ${\bold {MA}} \ \& \ \neg \CH$ . We can also get a  model for this
case by an $\omega_{2}$-iteration of infinitely equal reals over a model
for $\CH$.
\item[$ G$] $\omega_{2}$-iteration with finite support of eventually
different reals (see \cite{M1}) over a model for $\CH$.
\item[$ G^\star$] $\omega_{1}$-iteration with finite support of eventually
different reals  over a model for $\neg \CH \ \& \ {\bold B}(c)$.
\item[$ H\!=\!H^\star$] $\boldsymbol\aleph_{2}$ Cohen reals over a
  model for $\CH$. This model is self-dual. 
\item[$ I\!=\!I^\star$] $\omega_{2}$-iteration of Mathias forcing over
  a model for $\CH$  \cite{M1}. This model is self dual. 
\end{enumerate}

\begin{tabular}{|c|c|c|c|c|c|c|c|c|c|}\hline
\multicolumn{3}{|c|}
{\makebox[0.35in]{\rule[-0.2in]{0in}{0.5in}
$w{\bold D}\  \& \ \neg {\bold D}$
\rule[-0.2in]{0in}{0.5in}}}&
\makebox[0.35in]{\rule[-0.2in]{0in}{0.5in}
{\small Add}
\rule[-0.2in]{0in}{0.5in}}&
\makebox[0.35in]{\rule[-0.2in]{0in}{0.5in}
${\bold T}$
\rule[-0.2in]{0in}{0.5in}}&
\makebox[0.35in]{\rule[-0.2in]{0in}{0.5in}
${\bold F}$
\rule[-0.2in]{0in}{0.5in}}&
\makebox[0.35in]{\rule[-0.2in]{0in}{0.5in}
${\bold F}$
\rule[-0.2in]{0in}{0.5in}}&
\makebox[0.35in]{\rule[-0.2in]{0in}{0.5in}
${\bold F}$
\rule[-0.2in]{0in}{0.5in}}&
\makebox[0.35in]{\rule[-0.2in]{0in}{0.5in}
${\bold F}$
\rule[-0.2in]{0in}{0.5in}}&
\makebox[0.35in]{\rule[-0.2in]{0in}{0.5in}
${\bold F}$
\rule[-0.2in]{0in}{0.5in}}\\
\cline{4-10}
\multicolumn{3}{|r|}
{\makebox[0.35in]{\rule[-0.2in]{0in}{0.5in}
{\small Category}
\rule[-0.2in]{0in}{0.5in}}}&
\makebox[0.35in]{\rule[-0.2in]{0in}{0.5in}
{\small Baire}
\rule[-0.2in]{0in}{0.5in}}&
\makebox[0.35in]{\rule[-0.2in]{0in}{0.5in}
${\bold T}$
\rule[-0.2in]{0in}{0.5in}}&
\makebox[0.35in]{\rule[-0.2in]{0in}{0.5in}
${\bold T}$
\rule[-0.2in]{0in}{0.5in}}&
\makebox[0.35in]{\rule[-0.2in]{0in}{0.5in}
${\bold T}$
\rule[-0.2in]{0in}{0.5in}}&
\makebox[0.35in]{\rule[-0.2in]{0in}{0.5in}
${\bold F}$
\rule[-0.2in]{0in}{0.5in}}&
\makebox[0.35in]{\rule[-0.2in]{0in}{0.5in}
${\bold F}$
\rule[-0.2in]{0in}{0.5in}}&
\makebox[0.35in]{\rule[-0.2in]{0in}{0.5in}
${\bold F}$
\rule[-0.2in]{0in}{0.5in}}\\ \cline{4-10}
\multicolumn{3}{|c|}
{\makebox[0.35in]{\rule[-0.2in]{0in}{0.5in}
{\small Measure}
\rule[-0.2in]{0in}{0.5in}}}&
\makebox[0.35in]{\rule[-0.2in]{0in}{0.5in}
{\small Unif}
\rule[-0.2in]{0in}{0.5in}}&
\makebox[0.35in]{\rule[-0.2in]{0in}{0.5in}
${\bold T}$
\rule[-0.2in]{0in}{0.5in}}&
\makebox[0.35in]{\rule[-0.2in]{0in}{0.5in}
${\bold T}$
\rule[-0.2in]{0in}{0.5in}}&
\makebox[0.35in]{\rule[-0.2in]{0in}{0.5in}
${\bold F}$
\rule[-0.2in]{0in}{0.5in}}&
\makebox[0.35in]{\rule[-0.2in]{0in}{0.5in}
${\bold T}$
\rule[-0.2in]{0in}{0.5in}}&
\makebox[0.35in]{\rule[-0.2in]{0in}{0.5in}
${\bold F}$
\rule[-0.2in]{0in}{0.5in}}&
\makebox[0.35in]{\rule[-0.2in]{0in}{0.5in}
${\bold F}$
\rule[-0.2in]{0in}{0.5in}}\\ \hline
\makebox[0.35in]{\rule[-0.2in]{0in}{0.5in}
{\small Add}
\rule[-0.2in]{0in}{0.5in}}&
\makebox[0.35in]{\rule[-0.2in]{0in}{0.5in}
{\small Baire}
\rule[-0.2in]{0in}{0.5in}}&
\makebox[0.35in]{\rule[-0.2in]{0in}{0.5in}
{\small Unif}
\rule[-0.2in]{0in}{0.5in}}&
\makebox[0.35in]{\rule[-0.2in]{0in}{0.5in}
{\small Cov}
\rule[-0.2in]{0in}{0.5in}}&
\makebox[0.35in]{\rule[-0.2in]{0in}{0.5in}
${\bold T}$
\rule[-0.2in]{0in}{0.5in}}&
\makebox[0.35in]{\rule[-0.2in]{0in}{0.5in}
${\bold T}$
\rule[-0.2in]{0in}{0.5in}}&
\makebox[0.35in]{\rule[-0.2in]{0in}{0.5in}
${\bold T}$
\rule[-0.2in]{0in}{0.5in}}&
\makebox[0.35in]{\rule[-0.2in]{0in}{0.5in}
${\bold T}$
\rule[-0.2in]{0in}{0.5in}}&
\makebox[0.35in]{\rule[-0.2in]{0in}{0.5in}
${\bold T}$
\rule[-0.2in]{0in}{0.5in}}&
\makebox[0.35in]{\rule[-0.2in]{0in}{0.5in}
${\bold F}$
\rule[-0.2in]{0in}{0.5in}}\\ \hline
\makebox[0.35in]{\rule[-0.2in]{0in}{0.5in}
${\bold T}$
\rule[-0.2in]{0in}{0.5in}}&
\makebox[0.35in]{\rule[-0.2in]{0in}{0.5in}
${\bold T}$
\rule[-0.2in]{0in}{0.5in}}&
\makebox[0.35in]{\rule[-0.2in]{0in}{0.5in}
${\bold T}$
\rule[-0.2in]{0in}{0.5in}}&
\makebox[0.35in]{\rule[-0.2in]{0in}{0.5in}
${\bold T}$
\rule[-0.2in]{0in}{0.5in}}&
\makebox[0.35in]{\rule[-0.2in]{0in}{0.5in}
\grayrule
\rule[-0.2in]{0in}{0.5in}}&
\makebox[0.35in]{\rule[-0.2in]{0in}{0.5in}
\grayrule
\rule[-0.2in]{0in}{0.5in}}&
\makebox[0.35in]{\rule[-0.2in]{0in}{0.5in}
\grayrule
\rule[-0.2in]{0in}{0.5in}}&
\makebox[0.35in]{\rule[-0.2in]{0in}{0.5in}
\grayrule
\rule[-0.2in]{0in}{0.5in}}&
\makebox[0.35in]{\rule[-0.2in]{0in}{0.5in}
\grayrule
\rule[-0.2in]{0in}{0.5in}}&
\makebox[0.35in]{\rule[-0.2in]{0in}{0.5in}
\grayrule
\rule[-0.2in]{0in}{0.5in}}\\ \hline
\makebox[0.35in]{\rule[-0.2in]{0in}{0.5in}
${\bold F}$
\rule[-0.2in]{0in}{0.5in}}&
\makebox[0.35in]{\rule[-0.2in]{0in}{0.5in}
${\bold T}$
\rule[-0.2in]{0in}{0.5in}}&
\makebox[0.35in]{\rule[-0.2in]{0in}{0.5in}
${\bold T}$
\rule[-0.2in]{0in}{0.5in}}&
\makebox[0.35in]{\rule[-0.2in]{0in}{0.5in}
${\bold T}$
\rule[-0.2in]{0in}{0.5in}}&
\makebox[0.35in]{\rule[-0.2in]{0in}{0.5in}
\grayrule
\rule[-0.2in]{0in}{0.5in}}&
\makebox[0.35in]{\rule[-0.2in]{0in}{0.5in}
$A$
\rule[-0.2in]{0in}{0.5in}}&
\makebox[0.35in]{\rule[-0.2in]{0in}{0.5in}
\grayrule
\rule[-0.2in]{0in}{0.5in}}&
\makebox[0.35in]{\rule[-0.2in]{0in}{0.5in}
$ B$
\rule[-0.2in]{0in}{0.5in}}&
\makebox[0.35in]{\rule[-0.2in]{0in}{0.5in}
\grayrule
\rule[-0.2in]{0in}{0.5in}}&
\makebox[0.35in]{\rule[-0.2in]{0in}{0.5in}
\grayrule
\rule[-0.2in]{0in}{0.5in}}\\ \hline
\makebox[0.35in]{\rule[-0.2in]{0in}{0.5in}
${\bold F}$
\rule[-0.2in]{0in}{0.5in}}&
\makebox[0.35in]{\rule[-0.2in]{0in}{0.5in}
${\bold T}$
\rule[-0.2in]{0in}{0.5in}}&
\makebox[0.35in]{\rule[-0.2in]{0in}{0.5in}
${\bold F}$
\rule[-0.2in]{0in}{0.5in}}&
\makebox[0.35in]{\rule[-0.2in]{0in}{0.5in}
${\bold T}$
\rule[-0.2in]{0in}{0.5in}}&
\makebox[0.35in]{\rule[-0.2in]{0in}{0.5in}
\grayrule
\rule[-0.2in]{0in}{0.5in}}&
\makebox[0.35in]{\rule[-0.2in]{0in}{0.5in}
\grayrule
\rule[-0.2in]{0in}{0.5in}}&
\makebox[0.35in]{\rule[-0.2in]{0in}{0.5in}
\grayrule
\rule[-0.2in]{0in}{0.5in}}&
\makebox[0.35in]{\rule[-0.2in]{0in}{0.5in}
$ C$
\rule[-0.2in]{0in}{0.5in}}&
\makebox[0.35in]{\rule[-0.2in]{0in}{0.5in}
\grayrule
\rule[-0.2in]{0in}{0.5in}}&
\makebox[0.35in]{\rule[-0.2in]{0in}{0.5in}
\grayrule
\rule[-0.2in]{0in}{0.5in}}\\ \hline
\makebox[0.35in]{\rule[-0.2in]{0in}{0.5in}
${\bold F}$
\rule[-0.2in]{0in}{0.5in}}&
\makebox[0.35in]{\rule[-0.2in]{0in}{0.5in}
${\bold F}$
\rule[-0.2in]{0in}{0.5in}}&
\makebox[0.35in]{\rule[-0.2in]{0in}{0.5in}
${\bold T}$
\rule[-0.2in]{0in}{0.5in}}&
\makebox[0.35in]{\rule[-0.2in]{0in}{0.5in}
${\bold T}$
\rule[-0.2in]{0in}{0.5in}}&
\makebox[0.35in]{\rule[-0.2in]{0in}{0.5in}
\grayrule
\rule[-0.2in]{0in}{0.5in}}&
\makebox[0.35in]{\rule[-0.2in]{0in}{0.5in}
$ D$
\rule[-0.2in]{0in}{0.5in}}&
\makebox[0.35in]{\rule[-0.2in]{0in}{0.5in}
$ E\!=\!E^\star$
\rule[-0.2in]{0in}{0.5in}}&
\makebox[0.35in]{\rule[-0.2in]{0in}{0.5in}
$ F\!=\!F^\star$
\rule[-0.2in]{0in}{0.5in}}&
\makebox[0.35in]{\rule[-0.2in]{0in}{0.5in}
$ D^\star$
\rule[-0.2in]{0in}{0.5in}}&
\makebox[0.35in]{\rule[-0.2in]{0in}{0.5in}
\grayrule
\rule[-0.2in]{0in}{0.5in}}\\ \hline
\makebox[0.35in]{\rule[-0.2in]{0in}{0.5in}
${\bold F}$
\rule[-0.2in]{0in}{0.5in}}&
\makebox[0.35in]{\rule[-0.2in]{0in}{0.5in}
${\bold F}$
\rule[-0.2in]{0in}{0.5in}}&
\makebox[0.35in]{\rule[-0.2in]{0in}{0.5in}
${\bold F}$
\rule[-0.2in]{0in}{0.5in}}&
\makebox[0.35in]{\rule[-0.2in]{0in}{0.5in}
${\bold T}$
\rule[-0.2in]{0in}{0.5in}}&
\makebox[0.35in]{\rule[-0.2in]{0in}{0.5in}
\grayrule
\rule[-0.2in]{0in}{0.5in}}&
\makebox[0.35in]{\rule[-0.2in]{0in}{0.5in}
\grayrule
\rule[-0.2in]{0in}{0.5in}}&
\makebox[0.35in]{\rule[-0.2in]{0in}{0.5in}
\grayrule
\rule[-0.2in]{0in}{0.5in}}&
\makebox[0.35in]{\rule[-0.2in]{0in}{0.5in}
$ B^\star$
\rule[-0.2in]{0in}{0.5in}}&
\makebox[0.35in]{\rule[-0.2in]{0in}{0.5in}
$A^\star$
\rule[-0.2in]{0in}{0.5in}}&
\makebox[0.35in]{\rule[-0.2in]{0in}{0.5in}
\grayrule
\rule[-0.2in]{0in}{0.5in}}\\ \hline
\makebox[0.35in]{\rule[-0.2in]{0in}{0.5in}
${\bold F}$
\rule[-0.2in]{0in}{0.5in}}&
\makebox[0.35in]{\rule[-0.2in]{0in}{0.5in}
${\bold F}$
\rule[-0.2in]{0in}{0.5in}}&
\makebox[0.35in]{\rule[-0.2in]{0in}{0.5in}
${\bold F}$
\rule[-0.2in]{0in}{0.5in}}&
\makebox[0.35in]{\rule[-0.2in]{0in}{0.5in}
${\bold F}$
\rule[-0.2in]{0in}{0.5in}}&
\makebox[0.35in]{\rule[-0.2in]{0in}{0.5in}
\grayrule
\rule[-0.2in]{0in}{0.5in}}&
\makebox[0.35in]{\rule[-0.2in]{0in}{0.5in}
\grayrule
\rule[-0.2in]{0in}{0.5in}}&
\makebox[0.35in]{\rule[-0.2in]{0in}{0.5in}
\grayrule
\rule[-0.2in]{0in}{0.5in}}&
\makebox[0.35in]{\rule[-0.2in]{0in}{0.5in}
\grayrule
\rule[-0.2in]{0in}{0.5in}}&
\makebox[0.35in]{\rule[-0.2in]{0in}{0.5in}
\grayrule
\rule[-0.2in]{0in}{0.5in}}&
\makebox[0.35in]{\rule[-0.2in]{0in}{0.5in}
\grayrule
\rule[-0.2in]{0in}{0.5in}}\\ \hline
\end{tabular} 

\begin{enumerate}
\item[$A$] $\omega_{2}$-iteration with finite support of random reals over
a model for $\CH$.
\item[$ A^\star$] $\omega_{1}$-iteration with finite support of random
  reals over 
a model for $\neg \CH \ \& {\bold D}$.
\item[$ B$] $\omega_{2}$-iteration with countable support of forcing from
\cite{Sh2} and
random reals over
a model for $\CH$ (see section 5).
\item[$ B^\star$] $\omega_{2}$-iteration with countable support
of rational perfect set forcing and forcing ${\bold Q}_{f,g}$ from \cite{Sh3}
over
a model for $\CH$ (see section 5). 
\item[$ C$] $\boldsymbol\aleph_{2}$ Cohen and then
$\boldsymbol\aleph_{2}$ random reals over a model for $\CH$. This model is self-dual.
\item[$ D$] $\omega_{2}$-iteration of eventually different reals over a
model for $\CH$. \cite{M1}
\item[$ D^\star$] $\omega_{1}$-iteration of eventually different reals
  over a model for $\neg \CH \ \& \ {\bold B}(c)$.
\item[$ E\!=\!E^\star$] $\boldsymbol\aleph_{2}$ Cohen reals over a
  model for $\CH$. This model 
is self-dual.
\item[$ F\!=\!F^\star$] $\omega_{2}$-iteration with countable support of
forcing $Q$ from \cite{BS} over a model for $\CH$. This model is self-dual.
\end{enumerate}

Models in the following two tables are dual to each other.

\begin{tabular}{|c|c|c|c|c|c|c|c|c|c|}\hline
\multicolumn{3}{|l|}
{\makebox[0.35in]{\rule[-0.2in]{0in}{0.5in}
${\bold D}$
\rule[-0.2in]{0in}{0.5in}}}&
\makebox[0.35in]{\rule[-0.2in]{0in}{0.5in}
{\small Add}
\rule[-0.2in]{0in}{0.5in}}&
\makebox[0.35in]{\rule[-0.2in]{0in}{0.5in}
${\bold T}$
\rule[-0.2in]{0in}{0.5in}}&
\makebox[0.35in]{\rule[-0.2in]{0in}{0.5in}
${\bold F}$
\rule[-0.2in]{0in}{0.5in}}&
\makebox[0.35in]{\rule[-0.2in]{0in}{0.5in}
${\bold F}$
\rule[-0.2in]{0in}{0.5in}}&
\makebox[0.35in]{\rule[-0.2in]{0in}{0.5in}
${\bold F}$
\rule[-0.2in]{0in}{0.5in}}&
\makebox[0.35in]{\rule[-0.2in]{0in}{0.5in}
${\bold F}$
\rule[-0.2in]{0in}{0.5in}}&
\makebox[0.35in]{\rule[-0.2in]{0in}{0.5in}
${\bold F}$
\rule[-0.2in]{0in}{0.5in}}\\
\cline{4-10}
\multicolumn{3}{|r|}
{\makebox[0.35in]{\rule[-0.2in]{0in}{0.5in}
{\small Category}
\rule[-0.2in]{0in}{0.5in}}}&
\makebox[0.35in]{\rule[-0.2in]{0in}{0.5in}
{\small Baire}
\rule[-0.2in]{0in}{0.5in}}&
\makebox[0.35in]{\rule[-0.2in]{0in}{0.5in}
${\bold T}$
\rule[-0.2in]{0in}{0.5in}}&
\makebox[0.35in]{\rule[-0.2in]{0in}{0.5in}
${\bold T}$
\rule[-0.2in]{0in}{0.5in}}&
\makebox[0.35in]{\rule[-0.2in]{0in}{0.5in}
${\bold T}$
\rule[-0.2in]{0in}{0.5in}}&
\makebox[0.35in]{\rule[-0.2in]{0in}{0.5in}
${\bold F}$
\rule[-0.2in]{0in}{0.5in}}&
\makebox[0.35in]{\rule[-0.2in]{0in}{0.5in}
${\bold F}$
\rule[-0.2in]{0in}{0.5in}}&
\makebox[0.35in]{\rule[-0.2in]{0in}{0.5in}
${\bold F}$
\rule[-0.2in]{0in}{0.5in}}\\ \cline{4-10}
\multicolumn{3}{|c|}
{\makebox[0.35in]{\rule[-0.2in]{0in}{0.5in}
{\small Measure}
\rule[-0.2in]{0in}{0.5in}}}&
\makebox[0.35in]{\rule[-0.2in]{0in}{0.5in}
{\small Unif}
\rule[-0.2in]{0in}{0.5in}}&
\makebox[0.35in]{\rule[-0.2in]{0in}{0.5in}
${\bold T}$
\rule[-0.2in]{0in}{0.5in}}&
\makebox[0.35in]{\rule[-0.2in]{0in}{0.5in}
${\bold T}$
\rule[-0.2in]{0in}{0.5in}}&
\makebox[0.35in]{\rule[-0.2in]{0in}{0.5in}
${\bold F}$
\rule[-0.2in]{0in}{0.5in}}&
\makebox[0.35in]{\rule[-0.2in]{0in}{0.5in}
${\bold T}$
\rule[-0.2in]{0in}{0.5in}}&
\makebox[0.35in]{\rule[-0.2in]{0in}{0.5in}
${\bold F}$
\rule[-0.2in]{0in}{0.5in}}&
\makebox[0.35in]{\rule[-0.2in]{0in}{0.5in}
${\bold F}$
\rule[-0.2in]{0in}{0.5in}}\\ \hline
\makebox[0.35in]{\rule[-0.2in]{0in}{0.5in}
{\small Add}
\rule[-0.2in]{0in}{0.5in}}&
\makebox[0.35in]{\rule[-0.2in]{0in}{0.5in}
{\small Baire}
\rule[-0.2in]{0in}{0.5in}}&
\makebox[0.35in]{\rule[-0.2in]{0in}{0.5in}
{\small Unif}
\rule[-0.2in]{0in}{0.5in}}&
\makebox[0.35in]{\rule[-0.2in]{0in}{0.5in}
{\small Cov}
\rule[-0.2in]{0in}{0.5in}}&
\makebox[0.35in]{\rule[-0.2in]{0in}{0.5in}
${\bold T}$
\rule[-0.2in]{0in}{0.5in}}&
\makebox[0.35in]{\rule[-0.2in]{0in}{0.5in}
${\bold T}$
\rule[-0.2in]{0in}{0.5in}}&
\makebox[0.35in]{\rule[-0.2in]{0in}{0.5in}
${\bold T}$
\rule[-0.2in]{0in}{0.5in}}&
\makebox[0.35in]{\rule[-0.2in]{0in}{0.5in}
${\bold T}$
\rule[-0.2in]{0in}{0.5in}}&
\makebox[0.35in]{\rule[-0.2in]{0in}{0.5in}
${\bold T}$
\rule[-0.2in]{0in}{0.5in}}&
\makebox[0.35in]{\rule[-0.2in]{0in}{0.5in}
${\bold F}$
\rule[-0.2in]{0in}{0.5in}}\\ \hline
\makebox[0.35in]{\rule[-0.2in]{0in}{0.5in}
${\bold T}$
\rule[-0.2in]{0in}{0.5in}}&
\makebox[0.35in]{\rule[-0.2in]{0in}{0.5in}
${\bold T}$
\rule[-0.2in]{0in}{0.5in}}&
\makebox[0.35in]{\rule[-0.2in]{0in}{0.5in}
${\bold T}$
\rule[-0.2in]{0in}{0.5in}}&
\makebox[0.35in]{\rule[-0.2in]{0in}{0.5in}
${\bold T}$
\rule[-0.2in]{0in}{0.5in}}&
\makebox[0.35in]{\rule[-0.2in]{0in}{0.5in}
$ A$
\rule[-0.2in]{0in}{0.5in}}&
\makebox[0.35in]{\rule[-0.2in]{0in}{0.5in}
\grayrule
\rule[-0.2in]{0in}{0.5in}}&
\makebox[0.35in]{\rule[-0.2in]{0in}{0.5in}
\grayrule
\rule[-0.2in]{0in}{0.5in}}&
\makebox[0.35in]{\rule[-0.2in]{0in}{0.5in}
\grayrule
\rule[-0.2in]{0in}{0.5in}}&
\makebox[0.35in]{\rule[-0.2in]{0in}{0.5in}
\grayrule
\rule[-0.2in]{0in}{0.5in}}&
\makebox[0.35in]{\rule[-0.2in]{0in}{0.5in}
\grayrule
\rule[-0.2in]{0in}{0.5in}}\\ \hline
\makebox[0.35in]{\rule[-0.2in]{0in}{0.5in}
${\bold F}$
\rule[-0.2in]{0in}{0.5in}}&
\makebox[0.35in]{\rule[-0.2in]{0in}{0.5in}
${\bold T}$
\rule[-0.2in]{0in}{0.5in}}&
\makebox[0.35in]{\rule[-0.2in]{0in}{0.5in}
${\bold T}$
\rule[-0.2in]{0in}{0.5in}}&
\makebox[0.35in]{\rule[-0.2in]{0in}{0.5in}
${\bold T}$
\rule[-0.2in]{0in}{0.5in}}&
\makebox[0.35in]{\rule[-0.2in]{0in}{0.5in}
$ B$
\rule[-0.2in]{0in}{0.5in}}&
\makebox[0.35in]{\rule[-0.2in]{0in}{0.5in}
\grayrule
\rule[-0.2in]{0in}{0.5in}}&
\makebox[0.35in]{\rule[-0.2in]{0in}{0.5in}
\grayrule
\rule[-0.2in]{0in}{0.5in}}&
\makebox[0.35in]{\rule[-0.2in]{0in}{0.5in}
$ C$
\rule[-0.2in]{0in}{0.5in}}&
\makebox[0.35in]{\rule[-0.2in]{0in}{0.5in}
\grayrule
\rule[-0.2in]{0in}{0.5in}}&
\makebox[0.35in]{\rule[-0.2in]{0in}{0.5in}
\grayrule
\rule[-0.2in]{0in}{0.5in}}\\ \hline
\makebox[0.35in]{\rule[-0.2in]{0in}{0.5in}
${\bold F}$
\rule[-0.2in]{0in}{0.5in}}&
\makebox[0.35in]{\rule[-0.2in]{0in}{0.5in}
${\bold T}$
\rule[-0.2in]{0in}{0.5in}}&
\makebox[0.35in]{\rule[-0.2in]{0in}{0.5in}
${\bold F}$
\rule[-0.2in]{0in}{0.5in}}&
\makebox[0.35in]{\rule[-0.2in]{0in}{0.5in}
${\bold T}$
\rule[-0.2in]{0in}{0.5in}}&
\makebox[0.35in]{\rule[-0.2in]{0in}{0.5in}
\grayrule
\rule[-0.2in]{0in}{0.5in}}&
\makebox[0.35in]{\rule[-0.2in]{0in}{0.5in}
\grayrule
\rule[-0.2in]{0in}{0.5in}}&
\makebox[0.35in]{\rule[-0.2in]{0in}{0.5in}
\grayrule
\rule[-0.2in]{0in}{0.5in}}&
\makebox[0.35in]{\rule[-0.2in]{0in}{0.5in}
$ F$
\rule[-0.2in]{0in}{0.5in}}&
\makebox[0.35in]{\rule[-0.2in]{0in}{0.5in}
\grayrule
\rule[-0.2in]{0in}{0.5in}}&
\makebox[0.35in]{\rule[-0.2in]{0in}{0.5in}
\grayrule
\rule[-0.2in]{0in}{0.5in}}\\ \hline
\makebox[0.35in]{\rule[-0.2in]{0in}{0.5in}
${\bold F}$
\rule[-0.2in]{0in}{0.5in}}&
\makebox[0.35in]{\rule[-0.2in]{0in}{0.5in}
${\bold F}$
\rule[-0.2in]{0in}{0.5in}}&
\makebox[0.35in]{\rule[-0.2in]{0in}{0.5in}
${\bold T}$
\rule[-0.2in]{0in}{0.5in}}&
\makebox[0.35in]{\rule[-0.2in]{0in}{0.5in}
${\bold T}$
\rule[-0.2in]{0in}{0.5in}}&
\makebox[0.35in]{\rule[-0.2in]{0in}{0.5in}
$ E$
\rule[-0.2in]{0in}{0.5in}}&
\makebox[0.35in]{\rule[-0.2in]{0in}{0.5in}
\grayrule
\rule[-0.2in]{0in}{0.5in}}&
\makebox[0.35in]{\rule[-0.2in]{0in}{0.5in}
\grayrule
\rule[-0.2in]{0in}{0.5in}}&
\makebox[0.35in]{\rule[-0.2in]{0in}{0.5in}
$ D$
\rule[-0.2in]{0in}{0.5in}}&
\makebox[0.35in]{\rule[-0.2in]{0in}{0.5in}
\grayrule
\rule[-0.2in]{0in}{0.5in}}&
\makebox[0.35in]{\rule[-0.2in]{0in}{0.5in}
\grayrule
\rule[-0.2in]{0in}{0.5in}}\\ \hline
\makebox[0.35in]{\rule[-0.2in]{0in}{0.5in}
${\bold F}$
\rule[-0.2in]{0in}{0.5in}}&
\makebox[0.35in]{\rule[-0.2in]{0in}{0.5in}
${\bold F}$
\rule[-0.2in]{0in}{0.5in}}&
\makebox[0.35in]{\rule[-0.2in]{0in}{0.5in}
${\bold F}$
\rule[-0.2in]{0in}{0.5in}}&
\makebox[0.35in]{\rule[-0.2in]{0in}{0.5in}
${\bold T}$
\rule[-0.2in]{0in}{0.5in}}&
\makebox[0.35in]{\rule[-0.2in]{0in}{0.5in}
\grayrule
\rule[-0.2in]{0in}{0.5in}}&
\makebox[0.35in]{\rule[-0.2in]{0in}{0.5in}
\grayrule
\rule[-0.2in]{0in}{0.5in}}&
\makebox[0.35in]{\rule[-0.2in]{0in}{0.5in}
\grayrule
\rule[-0.2in]{0in}{0.5in}}&
\makebox[0.35in]{\rule[-0.2in]{0in}{0.5in}
$ G$
\rule[-0.2in]{0in}{0.5in}}&
\makebox[0.35in]{\rule[-0.2in]{0in}{0.5in}
\grayrule
\rule[-0.2in]{0in}{0.5in}}&
\makebox[0.35in]{\rule[-0.2in]{0in}{0.5in}
\grayrule
\rule[-0.2in]{0in}{0.5in}}\\ \hline
\makebox[0.35in]{\rule[-0.2in]{0in}{0.5in}
${\bold F}$
\rule[-0.2in]{0in}{0.5in}}&
\makebox[0.35in]{\rule[-0.2in]{0in}{0.5in}
${\bold F}$
\rule[-0.2in]{0in}{0.5in}}&
\makebox[0.35in]{\rule[-0.2in]{0in}{0.5in}
${\bold F}$
\rule[-0.2in]{0in}{0.5in}}&
\makebox[0.35in]{\rule[-0.2in]{0in}{0.5in}
${\bold F}$
\rule[-0.2in]{0in}{0.5in}}&
\makebox[0.35in]{\rule[-0.2in]{0in}{0.5in}
\grayrule
\rule[-0.2in]{0in}{0.5in}}&
\makebox[0.35in]{\rule[-0.2in]{0in}{0.5in}
\grayrule
\rule[-0.2in]{0in}{0.5in}}&
\makebox[0.35in]{\rule[-0.2in]{0in}{0.5in}
\grayrule
\rule[-0.2in]{0in}{0.5in}}&
\makebox[0.35in]{\rule[-0.2in]{0in}{0.5in}
\grayrule
\rule[-0.2in]{0in}{0.5in}}&
\makebox[0.35in]{\rule[-0.2in]{0in}{0.5in}
\grayrule
\rule[-0.2in]{0in}{0.5in}}&
\makebox[0.35in]{\rule[-0.2in]{0in}{0.5in}
\grayrule
\rule[-0.2in]{0in}{0.5in}}\\ \hline
\end{tabular}

\begin{tabular}{|c|c|c|c|c|c|c|c|c|c|}\hline
\multicolumn{3}{|l|}
{\makebox[0.35in]{\rule[-0.2in]{0in}{0.5in}
$\neg w{\bold D}$
\rule[-0.2in]{0in}{0.5in}}}&
\makebox[0.35in]{\rule[-0.2in]{0in}{0.5in}
{\small Add}
\rule[-0.2in]{0in}{0.5in}}&
\makebox[0.35in]{\rule[-0.2in]{0in}{0.5in}
${\bold T}$
\rule[-0.2in]{0in}{0.5in}}&
\makebox[0.35in]{\rule[-0.2in]{0in}{0.5in}
${\bold F}$
\rule[-0.2in]{0in}{0.5in}}&
\makebox[0.35in]{\rule[-0.2in]{0in}{0.5in}
${\bold F}$
\rule[-0.2in]{0in}{0.5in}}&
\makebox[0.35in]{\rule[-0.2in]{0in}{0.5in}
${\bold F}$
\rule[-0.2in]{0in}{0.5in}}&
\makebox[0.35in]{\rule[-0.2in]{0in}{0.5in}
${\bold F}$
\rule[-0.2in]{0in}{0.5in}}&
\makebox[0.35in]{\rule[-0.2in]{0in}{0.5in}
${\bold F}$
\rule[-0.2in]{0in}{0.5in}}\\
\cline{4-10}
\multicolumn{3}{|r|}
{\makebox[0.35in]{\rule[-0.2in]{0in}{0.5in}
{\small Category}
\rule[-0.2in]{0in}{0.5in}}}&
\makebox[0.35in]{\rule[-0.2in]{0in}{0.5in}
{\small Baire}
\rule[-0.2in]{0in}{0.5in}}&
\makebox[0.35in]{\rule[-0.2in]{0in}{0.5in}
${\bold T}$
\rule[-0.2in]{0in}{0.5in}}&
\makebox[0.35in]{\rule[-0.2in]{0in}{0.5in}
${\bold T}$
\rule[-0.2in]{0in}{0.5in}}&
\makebox[0.35in]{\rule[-0.2in]{0in}{0.5in}
${\bold T}$
\rule[-0.2in]{0in}{0.5in}}&
\makebox[0.35in]{\rule[-0.2in]{0in}{0.5in}
${\bold F}$
\rule[-0.2in]{0in}{0.5in}}&
\makebox[0.35in]{\rule[-0.2in]{0in}{0.5in}
${\bold F}$
\rule[-0.2in]{0in}{0.5in}}&
\makebox[0.35in]{\rule[-0.2in]{0in}{0.5in}
${\bold F}$
\rule[-0.2in]{0in}{0.5in}}\\ \cline{4-10}
\multicolumn{3}{|c|}
{\makebox[0.35in]{\rule[-0.2in]{0in}{0.5in}
{\small Measure}
\rule[-0.2in]{0in}{0.5in}}}&
\makebox[0.35in]{\rule[-0.2in]{0in}{0.5in}
{\small Unif}
\rule[-0.2in]{0in}{0.5in}}&
\makebox[0.35in]{\rule[-0.2in]{0in}{0.5in}
${\bold T}$
\rule[-0.2in]{0in}{0.5in}}&
\makebox[0.35in]{\rule[-0.2in]{0in}{0.5in}
${\bold T}$
\rule[-0.2in]{0in}{0.5in}}&
\makebox[0.35in]{\rule[-0.2in]{0in}{0.5in}
${\bold F}$
\rule[-0.2in]{0in}{0.5in}}&
\makebox[0.35in]{\rule[-0.2in]{0in}{0.5in}
${\bold T}$
\rule[-0.2in]{0in}{0.5in}}&
\makebox[0.35in]{\rule[-0.2in]{0in}{0.5in}
${\bold F}$
\rule[-0.2in]{0in}{0.5in}}&
\makebox[0.35in]{\rule[-0.2in]{0in}{0.5in}
${\bold F}$
\rule[-0.2in]{0in}{0.5in}}\\ \hline
\makebox[0.35in]{\rule[-0.2in]{0in}{0.5in}
{\small Add}
\rule[-0.2in]{0in}{0.5in}}&
\makebox[0.35in]{\rule[-0.2in]{0in}{0.5in}
{\small Baire}
\rule[-0.2in]{0in}{0.5in}}&
\makebox[0.35in]{\rule[-0.2in]{0in}{0.5in}
{\small Unif}
\rule[-0.2in]{0in}{0.5in}}&
\makebox[0.35in]{\rule[-0.2in]{0in}{0.5in}
{\small Cov}
\rule[-0.2in]{0in}{0.5in}}&
\makebox[0.35in]{\rule[-0.2in]{0in}{0.5in}
${\bold T}$
\rule[-0.2in]{0in}{0.5in}}&
\makebox[0.35in]{\rule[-0.2in]{0in}{0.5in}
${\bold T}$
\rule[-0.2in]{0in}{0.5in}}&
\makebox[0.35in]{\rule[-0.2in]{0in}{0.5in}
${\bold T}$
\rule[-0.2in]{0in}{0.5in}}&
\makebox[0.35in]{\rule[-0.2in]{0in}{0.5in}
${\bold T}$
\rule[-0.2in]{0in}{0.5in}}&
\makebox[0.35in]{\rule[-0.2in]{0in}{0.5in}
${\bold T}$
\rule[-0.2in]{0in}{0.5in}}&
\makebox[0.35in]{\rule[-0.2in]{0in}{0.5in}
${\bold F}$
\rule[-0.2in]{0in}{0.5in}}\\ \hline
\makebox[0.35in]{\rule[-0.2in]{0in}{0.5in}
${\bold T}$
\rule[-0.2in]{0in}{0.5in}}&
\makebox[0.35in]{\rule[-0.2in]{0in}{0.5in}
${\bold T}$
\rule[-0.2in]{0in}{0.5in}}&
\makebox[0.35in]{\rule[-0.2in]{0in}{0.5in}
${\bold T}$
\rule[-0.2in]{0in}{0.5in}}&
\makebox[0.35in]{\rule[-0.2in]{0in}{0.5in}
${\bold T}$
\rule[-0.2in]{0in}{0.5in}}&
\makebox[0.35in]{\rule[-0.2in]{0in}{0.5in}
\grayrule
\rule[-0.2in]{0in}{0.5in}}&
\makebox[0.35in]{\rule[-0.2in]{0in}{0.5in}
\grayrule
\rule[-0.2in]{0in}{0.5in}}&
\makebox[0.35in]{\rule[-0.2in]{0in}{0.5in}
\grayrule
\rule[-0.2in]{0in}{0.5in}}&
\makebox[0.35in]{\rule[-0.2in]{0in}{0.5in}
\grayrule
\rule[-0.2in]{0in}{0.5in}}&
\makebox[0.35in]{\rule[-0.2in]{0in}{0.5in}
\grayrule
\rule[-0.2in]{0in}{0.5in}}&
\makebox[0.35in]{\rule[-0.2in]{0in}{0.5in}
\grayrule
\rule[-0.2in]{0in}{0.5in}}\\ \hline
\makebox[0.35in]{\rule[-0.2in]{0in}{0.5in}
${\bold F}$
\rule[-0.2in]{0in}{0.5in}}&
\makebox[0.35in]{\rule[-0.2in]{0in}{0.5in}
${\bold T}$
\rule[-0.2in]{0in}{0.5in}}&
\makebox[0.35in]{\rule[-0.2in]{0in}{0.5in}
${\bold T}$
\rule[-0.2in]{0in}{0.5in}}&
\makebox[0.35in]{\rule[-0.2in]{0in}{0.5in}
${\bold T}$
\rule[-0.2in]{0in}{0.5in}}&
\makebox[0.35in]{\rule[-0.2in]{0in}{0.5in}
\grayrule
\rule[-0.2in]{0in}{0.5in}}&
\makebox[0.35in]{\rule[-0.2in]{0in}{0.5in}
\grayrule
\rule[-0.2in]{0in}{0.5in}}&
\makebox[0.35in]{\rule[-0.2in]{0in}{0.5in}
\grayrule
\rule[-0.2in]{0in}{0.5in}}&
\makebox[0.35in]{\rule[-0.2in]{0in}{0.5in}
$ G^\star$
\rule[-0.2in]{0in}{0.5in}}&
\makebox[0.35in]{\rule[-0.2in]{0in}{0.5in}
\grayrule
\rule[-0.2in]{0in}{0.5in}}&
\makebox[0.35in]{\rule[-0.2in]{0in}{0.5in}
\grayrule
\rule[-0.2in]{0in}{0.5in}}\\ \hline
\makebox[0.35in]{\rule[-0.2in]{0in}{0.5in}
${\bold F}$
\rule[-0.2in]{0in}{0.5in}}&
\makebox[0.35in]{\rule[-0.2in]{0in}{0.5in}
${\bold T}$
\rule[-0.2in]{0in}{0.5in}}&
\makebox[0.35in]{\rule[-0.2in]{0in}{0.5in}
${\bold F}$
\rule[-0.2in]{0in}{0.5in}}&
\makebox[0.35in]{\rule[-0.2in]{0in}{0.5in}
${\bold T}$
\rule[-0.2in]{0in}{0.5in}}&
\makebox[0.35in]{\rule[-0.2in]{0in}{0.5in}
\grayrule
\rule[-0.2in]{0in}{0.5in}}&
\makebox[0.35in]{\rule[-0.2in]{0in}{0.5in}
\grayrule
\rule[-0.2in]{0in}{0.5in}}&
\makebox[0.35in]{\rule[-0.2in]{0in}{0.5in}
\grayrule
\rule[-0.2in]{0in}{0.5in}}&
\makebox[0.35in]{\rule[-0.2in]{0in}{0.5in}
$ F^\star$
\rule[-0.2in]{0in}{0.5in}}&
\makebox[0.35in]{\rule[-0.2in]{0in}{0.5in}
\grayrule
\rule[-0.2in]{0in}{0.5in}}&
\makebox[0.35in]{\rule[-0.2in]{0in}{0.5in}
\grayrule
\rule[-0.2in]{0in}{0.5in}}\\ \hline
\makebox[0.35in]{\rule[-0.2in]{0in}{0.5in}
${\bold F}$
\rule[-0.2in]{0in}{0.5in}}&
\makebox[0.35in]{\rule[-0.2in]{0in}{0.5in}
${\bold F}$
\rule[-0.2in]{0in}{0.5in}}&
\makebox[0.35in]{\rule[-0.2in]{0in}{0.5in}
${\bold T}$
\rule[-0.2in]{0in}{0.5in}}&
\makebox[0.35in]{\rule[-0.2in]{0in}{0.5in}
${\bold T}$
\rule[-0.2in]{0in}{0.5in}}&
\makebox[0.35in]{\rule[-0.2in]{0in}{0.5in}
\grayrule
\rule[-0.2in]{0in}{0.5in}}&
\makebox[0.35in]{\rule[-0.2in]{0in}{0.5in}
\grayrule
\rule[-0.2in]{0in}{0.5in}}&
\makebox[0.35in]{\rule[-0.2in]{0in}{0.5in}
\grayrule
\rule[-0.2in]{0in}{0.5in}}&
\makebox[0.35in]{\rule[-0.2in]{0in}{0.5in}
$ D^\star$
\rule[-0.2in]{0in}{0.5in}}&
\makebox[0.35in]{\rule[-0.2in]{0in}{0.5in}
\grayrule
\rule[-0.2in]{0in}{0.5in}}&
\makebox[0.35in]{\rule[-0.2in]{0in}{0.5in}
$ E^\star$
\rule[-0.2in]{0in}{0.5in}}\\ \hline
\makebox[0.35in]{\rule[-0.2in]{0in}{0.5in}
${\bold F}$
\rule[-0.2in]{0in}{0.5in}}&
\makebox[0.35in]{\rule[-0.2in]{0in}{0.5in}
${\bold F}$
\rule[-0.2in]{0in}{0.5in}}&
\makebox[0.35in]{\rule[-0.2in]{0in}{0.5in}
${\bold F}$
\rule[-0.2in]{0in}{0.5in}}&
\makebox[0.35in]{\rule[-0.2in]{0in}{0.5in}
${\bold T}$
\rule[-0.2in]{0in}{0.5in}}&
\makebox[0.35in]{\rule[-0.2in]{0in}{0.5in}
\grayrule
\rule[-0.2in]{0in}{0.5in}}&
\makebox[0.35in]{\rule[-0.2in]{0in}{0.5in}
\grayrule
\rule[-0.2in]{0in}{0.5in}}&
\makebox[0.35in]{\rule[-0.2in]{0in}{0.5in}
\grayrule
\rule[-0.2in]{0in}{0.5in}}&
\makebox[0.35in]{\rule[-0.2in]{0in}{0.5in}
$ C^\star$
\rule[-0.2in]{0in}{0.5in}}&
\makebox[0.35in]{\rule[-0.2in]{0in}{0.5in}
\grayrule
\rule[-0.2in]{0in}{0.5in}}&
\makebox[0.35in]{\rule[-0.2in]{0in}{0.5in}
$ B^\star$
\rule[-0.2in]{0in}{0.5in}}\\ \hline
\makebox[0.35in]{\rule[-0.2in]{0in}{0.5in}
${\bold F}$
\rule[-0.2in]{0in}{0.5in}}&
\makebox[0.35in]{\rule[-0.2in]{0in}{0.5in}
${\bold F}$
\rule[-0.2in]{0in}{0.5in}}&
\makebox[0.35in]{\rule[-0.2in]{0in}{0.5in}
${\bold F}$
\rule[-0.2in]{0in}{0.5in}}&
\makebox[0.35in]{\rule[-0.2in]{0in}{0.5in}
${\bold F}$
\rule[-0.2in]{0in}{0.5in}}&
\makebox[0.35in]{\rule[-0.2in]{0in}{0.5in}
\grayrule
\rule[-0.2in]{0in}{0.5in}}&
\makebox[0.35in]{\rule[-0.2in]{0in}{0.5in}
\grayrule
\rule[-0.2in]{0in}{0.5in}}&
\makebox[0.35in]{\rule[-0.2in]{0in}{0.5in}
\grayrule
\rule[-0.2in]{0in}{0.5in}}&
\makebox[0.35in]{\rule[-0.2in]{0in}{0.5in}
\grayrule
\rule[-0.2in]{0in}{0.5in}}&
\makebox[0.35in]{\rule[-0.2in]{0in}{0.5in}
\grayrule
\rule[-0.2in]{0in}{0.5in}}&
\makebox[0.35in]{\rule[-0.2in]{0in}{0.5in}
$ A^\star$
\rule[-0.2in]{0in}{0.5in}}\\ \hline
\end{tabular} 

\begin{enumerate}
\item[$ A$] $\omega_{2}$-iteration of amoeba reals over a model for
  $\CH$ or any model for ${\bold {MA}}$. 
\item[$ A^\star$] $\omega_2$-iteration of amoeba reals over a model for
$\neg \CH$. 
\item[$ B$] $\omega_{2}$-iteration of dominating and random reals over a
model for $\CH$. \cite{M1}
\item[$ B^\star$] $\omega_2$-iteration of dominating and random reals
  over a model for $\neg \CH \ \& {\bold B}(c)$.
\item[$ C$] $\omega_{2}$-iteration with countable support of Mathias and
random reals (see section 5).
\item[$ C^\star$] $\omega_{2}$-iteration with countable support of forcing
${\bold Q}_{f,g}$ from \cite{Sh3}
(see section 2 and 3).
\item[$ D$] $\omega_{2}$-iteration with countable support
of Mathias reals over a model for
$\CH$.
\item[$ D^\star$] $\omega_{2}$-iteration with countable support
of ${\bold Q}_{f,g}$ and infinitely
equal reals over a model for $\CH$. (section 2)
\item[$ E$] $\omega_{2}$-iteration of dominating reals over a model for
$\CH$. \cite{M1}
\item[$ E^\star$] $\omega_2$-iteration of dominating reals over a
  model for $\neg \CH \ \& \ {\bold {MA}}$ or $\omega_{2}$-iteration
  with countable support of eventually equal reals.
\item[$ F$] $\boldsymbol\aleph_{2}$ random reals over a model for
  ${\bold {MA}} \ \& \ 2^{\boldsymbol\aleph_{0}} = \boldsymbol\aleph_{2}$.
\item[$ F^\star$] $\boldsymbol\aleph_{2}$ random reals over a model for $\CH$.
\item[$ G$] $\omega_{2}$-iteration with countable support
of Laver reals over a model for
$\CH$.
\item[$ G^\star$] $\omega_{2}$-iteration with countable support
of infinitely equal and random  reals over a model for
$\CH$. \cite{JS}
\end{enumerate}

\section {Not adding unbounded reals}
Our first goal is to construct a model for
$\ZFCa \ \& \ \neg w{\bold D}\ \&\        {\bold U}(c)\ \&\  \neg{\bold B}(m)\ \& \
{\bold U}(m)$.

We start with the definition of the forcing which will be used
in this construction. This family of forcing notions
was defined in \cite{Sh3}.

\begin{definition}
Let $f \in \omega^{\omega}$ and
$g \in \omega^{\omega \times \omega}$ be two functions
such that
\begin{enumerate}
\item $f(n) > \prod_{j<n} f(j)$ for $ n \in \omega$,
\item $g(n,j+1) > f(n)^{2}\cdot g(n,j)$ for $n,j \in \omega$,
\item $\min\{j \in \omega: g(n,j) > f(n+1)\}
\stackrel{n \rightarrow \infty}{\longrightarrow} \infty$.
\end{enumerate}
Let
$$\Seq^{f} = \bigcup_{n \in \omega} \prod_{j < n} f(j) .$$
For a tree $T$ define
$T^{[s]} = \{t \in T : s \subset t$ or $t \supset s\}$,
$\suc_{T}(s) = \{t \in T : t \supset s, \ \lh(t) = \lh(s)+1\}$. If
$T = T^{[s]}$ for some $s \in T$ then $s$ is called a {\em stem} of
$T$.

Let ${\bold Q}_{f,g}$ be the following notion of forcing:
$T \in {\bold Q}_{f,g}$ iff
\begin{enumerate}
\item $T$ is a perfect subtree of $\Seq^{f}$,
\item there exists a function $h \in \omega^{\omega}$ diverging to
  infinity such that
$$\exists n \ \forall m \geq n \ \forall s \in T \cap \omega^{m} \
|\suc_{T}(s)| \geq g(m,h(m)).$$ 
\end{enumerate}
Elements of
${\bold Q}_{f,g}$ are ordered by $\subseteq$.
\end{definition}

Let
${\bold Q}'_{f,g} \subset
{\bold Q}_{f,g}$ be the set defined as follows:
$T \in
{\bold Q}'_{f,g}$ iff
there exists $s_{0} \in \Seq^{f}$ such that $T=T^{[s_{0}]}$ and
there exists an increasing function $h \in \omega^{\omega}$ such that
$$\forall m \geq \lh(s_{0}) \ \forall s \in T \cap \omega^{m-1} \
|\suc_{T}(s)| \geq g(m,h(m)).$$

Clearly
${\bold Q}'_{f,g}$ is dense in
${\bold Q}_{f,g}$ and therefore from now on we will work with conditions
in this form.
Notice that
\begin{lemma}\label{2.1}
${\bold V}^{{\bold Q}_{f,g}} \models$``
${\bold V} \cap \omega^{\omega}$
is meager
in $\omega^{\omega}$''.
\end{lemma}
\Proof
Notice that if $r$ is a ${\bold Q}_{f,g}$-generic real then by an
easy density argument we show that
$$\forall h \in {\bold V} \cap \omega^{\omega} \ \forall^{\infty}n \
h(n) \neq r(n)  .$$
Therefore
${\bold V} \cap \omega^{\omega} \subset
\{h \in\omega^{\omega}: \forall^{\infty}n \
h(n) \neq r(n) \}$ which is a meager set.~$\QED$

\begin{definition}
We say that notion of forcing ${\bold P}$ is
$\omega^{\omega}$-bounding if
$$\forall \sigma \in {\bold V}^{{\bold P}} \cap \omega^{\omega} 
 \ \forces \exists r \in {\bold V} \cap \omega^{\omega} \forall n \
\sigma(n) \leq r(n)  .$$
\end{definition}

The following theorem was proved in \cite{Sh3}, we prove it here for
completeness; 
\begin{theorem} ${\bold Q}_{f,g}$ is $\omega^{\omega}$-bounding.
\end{theorem}
\Proof We will need the following

\begin{definition}
For $T, T' \in {\bold Q}_{f,g}$ and $\widehat{k} \in  \omega$
define
$T \geq_{\hat{k}} T'$ if
\begin{enumerate}
\item $T \geq T'$,
\item $\forall s \in T \ \suc_{T}(s) \neq \suc_{T'}(s) \rightarrow
|\suc_{T}(s)| \geq g(\lh(s), \widehat{k})$.
\end{enumerate}
\end{definition}

\begin{claim}
Suppose that $\{T^{n} : n \in \omega\}$ is a sequence of elements
of ${\bold Q}_{f,g}$ such that $T^{n+1} \geq_{k_{n}} T^{n}$ for
$n \in \omega$ where $\{k_{n}: n \in \omega\}$ is an increasing sequence
of natural numbers.
Then there exists $T \in {\bold Q}_{f,g}$ such that
$T \geq_{k_{n}} T^{n}$ for $n \in \omega$.
\end{claim}
\Proof
For $n \in \omega$ define
$$u_{n} = \min\{j \in \omega: \forall k \geq j \ \forall s \in T^{n} \cap
\omega^{k} \ |\suc_{T_{n}}(s)| \geq g(k,k_{n})\} .$$
Let $T = \bigcup_{n \in \omega} T^{n} \rest u_{n}$.
Function $h(m)=k_{n-1}$ for $m \in [u_{n-1},u_n)$ witnesses
that
$T \in {\bold Q}_{f,g}$.~$\QED$

\begin{lemma}
Let
$T \in {\bold Q}_{f,g}$ and $\tau$ be such that $T \forces \tau \in \omega$.
Suppose that $\widehat{k} \in \omega$. Then there exists $\widehat{T} \geq_{\hat{k}} T$
and $n \in \omega$ such that
$$\forall s \in \widehat{T} \cap \omega^{n} \ \exists a_{s} \in \omega \
\widehat{T}^{[s]} \forces \tau = a_{s} .$$
\end{lemma}
\Proof
Let $S \subseteq T$ be the set of all $t \in  T$ such that
$T^{[t]}$ satisfies the lemma. In other words
$$S = \{ t  \in T: \exists n_{t} \in \omega \
\exists \widehat{T} \geq_{\hat{k}} T^{[t]} \
\forall s \in \widehat{T} \cap \omega^{n_{t}} \ \exists a_{s} \in \omega
\ \widehat{T}^{[s]} \forces \tau =a_{s}\} .$$
We want to show that stem of $T$ belongs to $S$.
Notice that if $s \not \in S$ then
$$|\suc_{T}(s) \cap S| \leq g(\lh(s),\widehat{k}).$$

Suppose that stem of $T$ does not belong to $S$ and by induction on levels
build a tree $\widehat{S} \geq_{\hat{k}} T$ such that for
$s \in \widehat{S}$,
$$\suc_{\hat{S}}(s) = \left\{ \begin{array}{ll} \suc_{T}(s) & \hbox{if }
|\suc_{T}(s) \cap S| \leq g(\lh(s),\widehat{k}) \\
\suc_{T}(s) - \suc_{S}(s) & \hbox{otherwise}
\end{array} \right. .$$
Clearly $\widehat{S} \in {\bold Q}_{f,g}$ since $g(\lh(s),m) -
g(\lh(s),\widehat{k}) \geq g(\lh(s),m-\widehat{k})$ for all $s$ and $m >
\widehat{k}$. 

Find $\widehat{S}_{1} \geq \widehat{S} $ and $\widehat{n} \in \omega$
such that 
$\widehat{S}_{1} \forces \tau=\widehat{n}$. Now get $t \in T$ and $\widehat{S}_{2} \geq
\widehat{S}_{1}$ such that $\widehat{S}_{2} \geq_{\hat{k}} T^{[t]}$.
But that contradicts  the definition of the condition
$\widehat{S}$.~$\QED$ 

We finish the proof of the theorem. Suppose that $T \forces \sigma \in
\omega^{\omega}$.
Build by induction sequences $\{T_{n}: n \in \omega\}$ and
$\{k_{n} : n \in \omega\}$ such that for $n \in \omega$,
\begin{enumerate}
\item $T_{n+1} \geq_{k_{n}} T_{n}$,
\item $\forall s \in T_{n+1} \cap \omega^{k_{n}} \ \exists a_{s} \in \omega
\ T_{n+1}^{[s]} \forces \sigma(n) = a_{s}$.
\end{enumerate}
Let $T = \lim_{n \rightarrow \infty} T_{n}$ and let
$r(n) = \max\{a_{s} : s \in T \cap \omega^{k_{n}}\}$ for $n \in \omega$.
Then
$$T \forces \forall n \in \omega \ \sigma(n) \leq r(n) $$
which finishes the proof.~$\QED$

Notice that in fact we proved that
\begin{lemma}\label{crucial}
If $T \forces \sigma \in \omega^{\omega}$ then there exists a sequence
$\{k_{n} : n \in \omega\} $ and a tree $\widehat{T} \geq T$ such that
$$\forall s \in \widehat{T} \cap \omega^{k_{n}} \ \exists a_{s} \in \omega \
\widehat{T}^{[s]} \forces \sigma(n) = a_{s}  .\ \QED$$
\end{lemma}

Our next goal is to show that forcing with ${\bold Q}_{f,g}$ does
not add random reals.
We will need the following

\begin{definition}
Let $f \in \omega^{\omega}$ and let $X_{f} = \prod_{n=0}^{\infty} f(n)$.
Define $S_{f}$ as follows:
$T \in S_{f}$ if
$T$ is a perfect subtree of $\Seq^{f}$ and
$$\lim_{n \rightarrow \infty}
\frac{|T \cap \omega^{n}|}
{\prod_{m=1}^{n-1} f(m)}
= 0 .$$
Notion of forcing ${\bold Q}$ is called $f$-bounding
if
$$\forall \sigma \in X_{f} \cap {\bold V}^{{\bold Q}}\
\exists T \in S_{f}\cap {\bold V} \
\forall n \ \sigma \rest n \in T  .$$
\end{definition}

\begin{theorem} Let ${\bold P}$ be a notion of forcing. We have the following
\begin{enumerate}
\item
If ${\bold P}$ is an $f$-bounding notion of forcing then
${\bold P}$ does not add random reals.
\item If ${\bold P}$ is $\omega^{\omega}$-bounding and ${\bold P}$
does not random reals then ${\bold P}$ is $f$-bounding for every
$f \in \omega^{\omega}$.
\end{enumerate}
\end{theorem}
\Proof
Define a measure $\mu$ on $X_{f}$ as a product of equally distributed,
normalized measures on $f(n)$.

(1)
Every element of $S_{f}$ corresponds to a closed,
measure zero subset of $X_{f}$.
This finishes the proof as $X_{f}$ is isomorphic to the Cantor space
with standard measure.

(2) Suppose that $\forces \sigma \in X_{f}$. Since we assume that
${\bold P}$ does not add random reals we can find a null $G_{\delta}$ subset
$H \in {\bold V}$ of $X_{f}$ such that $\forces \sigma \in H$.

\begin{claim} Suppose that $H \subseteq X_{f}$. Then
$\mu(H)=0$  iff there exists a sequence
$\{J_{n} \subseteq \Seq^{f} \cap \omega^{n} : n \in \omega \}$ such that
\begin{enumerate}
\item $H \subseteq \{x \in X_{f}: \exists^{\infty}n \ x \rest n
\in J_{n} \} \ ,$
\item $\sum_{n=0}^{\infty} \mu(\{x \in X_{f}:
x \rest n \in  J_{n}\}) < \infty$ .
\end{enumerate}
\end{claim}
\Proof
$(\leftarrow)$ This implication is an
immediate consequence of Borel-Cantelli lemma.

$(\rightarrow)$ Since $\mu(H)=0$ there are
open sets $\{G_{n} :n \in \omega\}$ covering $H$ such
that $\mu(G_{n}) < \frac{1}{2^{n}}$
for $n \in \omega$.
Write each $G_{n}$ as a union of disjoint basic sets i.e.
$$G_{n} = \bigcup_{m \in \omega}[s^{n}_{m}] \   for \  n \in \omega .$$
Let 
$J_{n} = \{s \in \Seq^{f} \cap \omega^{n} : s=s_{k}^{l}$ for some $ k,l \in
\omega \}$
for
$n \in \omega$.
Verification of (1) and (2) is straightforward.~$\QED$

Let $\{J_{n} : n \in \omega\}$ be a sequence obtained by applying the above
to the set $H$. In particular
$\{ n \in \omega : \sigma \rest n \in J_{n} \}$
is infinite. Using the fact that forcing
${\bold P}$ is $\omega^{\omega}$-bounding find a function $h \in \omega^{\omega}$
such that
$\forall n \ \exists m \in [h(n),h(n+1)) \ \sigma
\rest m \in J_{m}  .$
Let
$$C = \bigcap_{n \in \omega} \bigcup_{m=h(n)}^{h(n+1)} \bigcup_{s \in J_{m}}
[s]  .$$ 
It is easy to see that $C$ is a closed set and that
$\forces \sigma \in C$. As $C$ is a closed set $C$ is a set of branches of
some tree $T$. This tree has required properties.~$\QED$

The following theorem was proved in \cite{Sh3}, we prove it here for
completness. 
\begin{theorem}
Forcing ${\bold Q}_{f,g}$ is $f$-bounding.
\end{theorem}
\Proof
We start with the following
\begin{lemma}\label{2.9}
If $\widetilde{T} \forces \forall n \ \sigma(n) \leq f(n)$ then there exists
tree $\widehat{T} \geq \widetilde{T}$ such that
$$\forall s \in \widehat{T} \cap \omega^{n} \ \exists a_{s} \leq f(n) \
\widehat{T}^{[s]} \forces \sigma(n) = a_{s} .$$
\end{lemma}
\Proof
By applying \ref{crucial} we get a tree $T \geq \widetilde{T}$ and a sequence
$\{k_{n} : n \in \omega\}$ such that
$$\forall s \in T \cap \omega^{k_{n}} \ \exists a_{s} \in \omega \
T^{[s]} \forces \sigma(n) = a_{s} .$$
Without loss of generality we can assume that
$k_{n} \geq n$ for all $n \in \omega$.
Suppose that function $h \in \omega^{\omega}$ witnesses that
$T \in {\bold Q}_{f,g}$. In other words $|\suc_{T}(s)| \geq
g(\lh(s),h(\lh(s)))$ 
for $s \in T$.

Build by induction a family of trees $\{T_{n,l} : n \in \omega, n \leq
l       \leq 
k_{n}\}$ such that
\begin{enumerate}
\item $T_{n,l} \geq T_{n,l'}$
for $l \leq l', \ n \in \omega$,
\item $T_{n,l} \rest l = T_{n,l'} \rest l$
for $l \leq l', \ n \in \omega$,
\item $T_{n,l} \geq T_{m,l'} $
for $n < m$ and all $l,l'\in \omega$,
\item $T_{n,l} \rest n = T_{m,l'} \rest n$
for $n < m$ and all $l,l'$,
\item $\forall n \
\forall s \in T_{n,l} \cap \omega^{l} \ \exists a_{s} \leq f(n) \
T_{n,l}^{[s]} \forces \sigma(n)=a_{s}$,
\item $\forall n \
\forall s \in T_{n,n} \cap \omega^{\leq n} \ |\suc_{T_{n,l}}(s)| \geq
g(\lh(s),h(\lh(s))-1)$.
\end{enumerate}
It is clear that $$\widehat{T} = 
\lim_{n \rightarrow \infty} T_{n,n}$$
has the required properties and the function $h'(n) = h(n)-1$
witnesses that 
$\widehat{T} \in {\bold Q}_{f,g}$.

Suppose that the tree $T_{n,n}$ is given for some $n \in \omega$.
Trees
$T_{n+1,k_{n}} \geq
T_{n+1,k_{n}-1} \geq
\ldots \geq
T_{n+1,n+1}$ are constucted by induction as follows:

Let $T_{n+1,k_{n}} = T_{n,n}$ and suppose that
$T_{n+1,l}$ is given.
Tree
$T_{n,l-1}$
will be defined in the following way:
$T_{n,l-1}\rest l-1 =
T_{n,l}\rest l-1$ and for each $t \in T_{n,l} \cap \omega^{l-1}$
we will specify which of the immediate successors of $t$ belong to
$T_{n,l-1}$.

Take $t \in T_{n+1,l} \cap \omega^{l-1}$ and let
$s \in \suc_{T_{n+1,l}}(t)$.
By (5) there exists $a_{s} \leq f(n) $ such that
$T_{n+1,l}^{[s]} \forces \sigma(n) = a_{s}$.
That defines a partition of the set
$\suc_{T_{n+1,l}}(t)$ into $f(n)$ many pieces.
Let the set of immediate successors of $t$ in
$T_{n+1,l-1}$ be the largest piece in this partition.

Notice that for $t \in T \cap \omega^{n}$ the set $\suc_{T}(t)$ will
be altered at most $n$ times and each time its size will decrease
by a factor $f(i)$ for $i \leq n$. Therefore
$$|\suc_{T_{n,n}}(t)| > \frac{g(n,h(n))}{\prod_{i\leq n} f(i)}
\geq g(n,h(n)-1) .$$
This verifies (6) and finishes the proof of the lemma.~$\QED$

Now we can prove the theorem.
Let $\sigma$ be a
${\bold Q}_{f,g}$-name
such that $\widetilde{T} \forces \forall n \ \sigma(n) \leq f(n)$
for some $\widetilde{T} \in
{\bold Q}_{f,g}$.

Let $\widehat{T} \geq \widetilde{T}$ be the condition as in the lemma
above. 
The tree $T'$ we are looking for will be defined as follows:
$$s \in T'  \mbox{ iff }  \exists t \in \widehat{T} \
\widehat{T}^{[t]} \forces \sigma \rest \lh(s) = s   .$$
By trimming $\widehat{T}$ some more we can see that
$$\frac{|T' \cap \omega^{n}|}
{\prod_{m=1}^{n} f(m)} \leq
\frac{|\widehat{T} \cap \omega^{n}|}
{\prod_{m=1}^{n} f(m)}
\stackrel{n \rightarrow \infty}
{\longrightarrow 0}.~\QED$$

To conclude this section we need some preservation theorems. We have
to show that a countable support iteration of $\omega^{\omega}$-bounding
forcings is $\omega^{\omega}$-bounding.
This has been proved for proper forcings (see \cite{Sh1}).
Here we present a much easier proof that works for a more limited class
of partial orderings.
Similarly we need to know that the iterations we use do not add
random reals. Unfortunately $f$-boundedness is not preserved by
a countable support iteration. We will prove it only for certain partial
orderings. For a general preservation theorem of a  slightly stronger
property called $(f,g)$-boundedness see \cite{Sh4}.

\begin{definition}
Let ${\bold P}$ be a notion of forcing satisfying axiom A (see
\cite{Bau}).   
We say that
${\bold P}$ has property ($\star$) if
for every $p \in {\bold P}, \widehat{n} \in \omega$ and a ${\bold P}$-name
$\tau$ for
a natural number there exists $N \in \omega$ and $q \geq_{\hat{n}} p$
such that $q \forces \tau < N$.
\end{definition}

It is easy to see that partial orderings having property ($\star$)
are $\omega^{\omega}$-bounding.
\begin{theorem}\label{bound}
Let $\{{\bold P}_{\xi}, \dot{{\bold Q}}_{\xi} : \xi < \alpha\}$ be a countable
support iteration of forcings that have the  property $(\star)$.
Then ${\bold P}_{\alpha} = \lim_{\xi < \alpha} {\bold P}_{\xi}$ is
$\omega^{\omega}$-bounding.
\end{theorem}
\Proof
For $p,q \in {\bold P}_{\alpha}, \ F \in [\alpha]^{<\omega}$ and
$\widehat{n} \in \omega$ write
$p \geq_{F,\hat{n}} q$ if
\begin{enumerate}
\item $p \geq q$,
\item $\forall \xi \in F \ p \rest \xi \forces p(\xi) \geq_{\hat{n}}
q(\xi)$.
\end{enumerate}

The proof of the theorem is based on the following general fact:
\begin{lemma}\label{bounding}
Suppose that $p \in {\bold P}_\alpha, \ F \in [\alpha]^{<\omega}$ and
$\widehat{n} \in \omega$
are given. Let $\tau$ be a ${\bold P}$-name for a natural number.
Then there exists $q \geq_{F,\hat{n}} p$ and $N \in \omega$
such that $q \forces \tau < N$.
\end{lemma}
\Proof
It will be proved by induction on $(|F|, \min F)$ over all possible
models. 
Suppose that $|F|=n+1$ and $\min F = \alpha_{0} < \alpha$.
By induction hypothesis in
${\bold V}^{P_{\alpha_{0}+1}}$
the lemma is
true for $F' = F - \{\alpha_{0}\}$.
Therefore there exists a ${\bold Q}_{\alpha_{0}}$ name $\sigma \in
{\bold V}^{P_{\alpha_0}}$
such that
$${\bold V}^{P_{\alpha_{0}}} \models \hbox{``}p(\alpha_{0}) \forces \exists
q \hbox{''} \geq_{F',\hat{n}} p \rest (\alpha_{0},\alpha) \ q
\hbox{''}  \forces
\tau < \sigma \hbox{''} .$$
Since ${\bold Q}_{\alpha_{0}}$ has property $(\star)$ in
${\bold V}^{P_{\alpha_{0}}}$ we can find $q' \geq_{\hat{n}}
p(\alpha_{0})$ and $N$ such that 
$${\bold V}^{P_{\alpha_{0}}} \models \hbox{``} q' \forces \sigma < N
\hbox{''} .$$ 
The last statement is forced by a condition $q_0 \in {\bold P}_{\alpha_0}$.
Let $q = q_0^\frown q^{' \frown} q''$. It is the condition we were
looking for.~$\QED$

Let $p_{0}$ be any element of ${\bold P}_{\alpha}$. Suppose that
$p_{0} \forces \sigma \in \omega^{\omega}$.
Using \ref{bounding} define by induction sequences
$\{p_{n} : n \in \omega\}, \{F_{n}: n \in \omega\}$ and
a function $r \in \omega^{\omega}$ such that
\begin{enumerate}
\item $p_{n+1} \geq_{F_{n},n} p_{n}$ for $ n \in \omega$,
\item $\forall \xi \in \supp(p_{n}) \ \exists j \in \omega \ \xi \in
  F_{j}$, 
\item $F_{n} \subset F_{n+1}$ for $n \in \omega$,
\item $p_{n+1} \forces \sigma(n) < r(n)$.
\end{enumerate}
Let $q $ be the limit of $\{p_{n}: n \in \omega\}$. Then
$q \forces \forall n \in \omega \ \sigma(n) < r(n)$.~$\QED$

Finally we can prove:
\begin{theorem}
$\Con(\ZFCa) \rightarrow \Con(\ZFCa \ \& \ \neg w{\bold D} \ \&\ {\bold U}(c) \ \&\
\neg{\bold B}(m) \ \&\ {\bold U}(m))$.
\end{theorem}
\Proof
The following notion of forcing was introduced in \cite{M1}:
let $f \in \omega^{\omega}$. Define

$p \in {\bold Q}_{f}$ iff
\begin{enumerate}
\item $p :\dom(p) \longrightarrow \omega$,
\item $\dom(p) \subset \omega$ and
$\omega- \dom(p)$
is infinite,
\item $\forall n \ p(n) \leq f(n)$.
\end{enumerate}
For $p,q \in {\bold Q}_{f}$ $p \geq q$ if $p \supseteq q$ and
for $n \in \omega \ p \geq_{n} q$ iff $p \geq q$ and the first $n$ elements
of
$\omega- \dom(p)$ and
$\omega- \dom(q)$ are the same.

The following fact is well known:
\begin{lemma}
Let ${\bold P}$ be a notion of forcing.
If ${\bold P}$ has the Laver property then ${\bold P}$ is
$f$-bounding for all  functions
$f \in \omega^{\omega}$.~$\QED$
\end{lemma}

\begin{lemma}\label{2.14} 
Let $f \in \omega^{\omega}$ be a strictly increasing function
such that $f(n) > 2^{n}$ for $n \in \omega$. Then
\begin{enumerate}
\item ${\bold V} \cap 2^{\omega}$ has measure zero in
${\bold V}^{{\bold Q}_{f}}$,
\item ${\bold Q}_{f}$ is $f$-bounding.
\end{enumerate}
\end{lemma}
\Proof
(1) It is enough to show that
$X_{f} \cap {\bold V}$
has measure zero in
${\bold V}^{{\bold Q}_{f}}$.
Notice that for $h \in X_{f}$ the set
$$H_{h} = \{x \in X_{f} : \exists^{\infty} n \ x(n) = h(n)\}$$
has measure zero.
It is easy to see that
$$X_{f} \cap {\bold V} \subset
H_{h_{G}} $$
where $h_{G}$ is a generic real.

(2)
Let $p_{0}$ be any element of ${\bold Q}_{f}$. Suppose that
$p_{0} \forces \sigma \in X_{f}$.
Define by induction sequences
$\{p_{n} : n \in \omega\}, \ \{k_{n} : n \in \omega\}$
and $\{J_{n} : n \in \omega\}$
such that
\begin{enumerate}
\item $J_{n} \subset \Seq^{f} \cap \omega^{k_{n}}$ for $n \in \omega$,
\item $p_{n+1} \geq_{n} p_{n}$ for $ n \in \omega$,
\item $p_{n+1} \forces \sigma \rest k_{n} \in J_{n}$ for $n \in \omega$,
\item $\displaystyle \frac{|J_{n}|}{\prod_{m=1}^{k_{n}} f(m)} \leq
  \displaystyle \frac{1}{n}$ for $n \in \omega$.
\end{enumerate}
Let $q \geq p_{0} $ be the limit of $\{p_{n}: n \in \omega\}$ and
$T = \bigcup_{n \in \omega} J_{n}$. By removing all nodes whose
ancestors are missing we can make sure that $T$ is a tree.
Then
$q$ forces that $\sigma$ is a branch through $T$ and by (4) $T$ has measure
zero.~$\QED$

Let  $\{{\bold P}_{\xi}, \dot{{\bold Q}}_{\xi} :
\xi < \boldsymbol\aleph_{2}\}$ be a countable support iteration such that

$\forces_{\xi}$ ``$\dot{{\bold Q}}_{\xi} \cong {\bold Q}_{f,g}$'' if $\xi$ is even

$\forces_{\xi}$ ``$\dot{{\bold Q}}_{\xi}\cong {\bold Q}_{f}$'' if $\xi$ is odd.

Let ${\bold P} ={\bold P}_{\boldsymbol\aleph_{2}}$. Then
${\bold V}^{\bold P} \ \models \
\neg w{\bold D}$
since ${\bold P}$ is $\omega^{\omega}$-bounding, 
and
${\bold V}^{\bold P} \models
{\bold U}(c) \ \& \ {\bold U}(m)$ by the properties of forcings
${\bold Q}_{f,g}$  and ${\bold Q}_f$ (note that ${\bold Q}_{f,g}$ has property
$(\star)$).
To finish the proof we need 
\begin{lemma}\label{preserve}
${\bold P}$ is $f$-bounding.
\end{lemma}
\Proof
For $p,q \in {\bold P}, \ F \in [\boldsymbol\aleph_{2}]^{<\omega}$ and
$\widehat{n} \in \omega$ denote
$p \geq_{F,\hat{n}} q$ if
\begin{enumerate}
\item $p \geq q$,
\item $\forall \xi \in F \ p \rest \xi \forces p(\xi) \geq_{\hat{n}}
q(\xi)$.
\end{enumerate}

Let $p_{0}$ be any element of ${\bold P}$. Suppose that
$p_{0} \forces \sigma \in X_{f}$.
Using the fact that both ${\bold Q}_{f,g}$  and  ${\bold Q}_{f}$ are
$f$-bounding and arguing as in the proofs of \ref{2.9} and \ref{2.14},
define by induction sequences
$\{p_{n} : n \in \omega\}, \{F_{n}: n \in \omega\},\ \{k_{n} : n \in \omega\}$
and $\{J_{n} : n \in \omega\}$
such that
\begin{enumerate}
\item $J_{n} \subset \Seq^{f} \cap \omega^{k_{n}}$ for $n \in \omega$,
\item $p_{n+1} \geq_{F_{n},n} p_{n}$ for $ n \in \omega$,
\item $\forall \xi \in \supp(p_{n}) \ \exists j \in \omega \ \xi \in F_{j}$,
\item $F_{n} \subset F_{n+1}$ for $n \in \omega$,
\item $p_{n+1} \forces \sigma \rest k_{n} \in J_{n}$ for $n \in \omega$,
\item $\frac{|J_{n}|}{\prod_{m=1}^{k_{n}} f(m)} \leq \frac{1}{n}$ for
$n \in \omega$.
\end{enumerate}
Let $q\geq p_{0} $ be the limit of $\{p_{n}: n \in \omega\}$ and
$T = \bigcup_{n \in \omega} J_{n}$. As before, by removing
non-splitting nodes we can assume that $T$ is a tree.
Then
$q$ forces that $\sigma$ is a branch through $T$ and by (6) $T$ has measure
zero.~$\QED$

Notice that \ref{preserve} can be proved in the same way for
many other forcings including perfect set forcing from section 5.

\section{Preserving ``old reals have outer measure 1''}

In this section we construct a model for
$\ZFCa \ \& \ \neg w{\bold D}\ \&\
{\bold U}(c)\ \&\ \neg{\bold U}(m)\ \&\ \neg{\bold B}(m)$.
It is obtained by $\omega_{2}$-iteration with countable support
of ${\bold Q}_{f,g}$.

The main problem is to verify that
$\neg {\bold U}(m)$ holds in that model.

We will use the following technique from \cite{JS}.

\begin{definition}
Let ${\bold P}$ be a notion of forcing. Define

$\star_{1}[{\bold P}]$
iff for every sufficiently large cardinal $\kappa$,
and for every countable elementary submodel $N \prec H(\kappa,\in)$, if
${\bold P} \in N$ and $\{\dot{I}_{n} : n \in \omega\} \in N$ is a
${\bold P}$-name for a sequence of rational intervals and
$\{p_{n}: n \in \omega\} \in N$ is a sequence of elements of ${\bold P}$
such that $p_{0} \forces \sum_{n=1}^{\infty} \mu(\dot{I}_{n})<\infty$ and
$p_{n} \forces \dot{I}_{n} = I_{n}$ for $n \in \omega$ then for
every random real $x$ over $N$, if $x \not \in \bigcup_{n \in \omega} I_{n}$
then there exists $q \geq p_{0}$ such that
\begin{enumerate}
\item $q$ is $(N,{\bold P})$-generic,
\item $q \forces x$ is random over $N[G]$ for every ${\bold P}$-generic filter
over $N$ containing $p_{0}$,
\item $q \forces x \not \in \bigcup_{n \in \omega} \dot{I}_{n}$.
\end{enumerate}

$\star_{2}[{\bold P}]$ iff
for every ${\bold P}$-name $\dot{A}$ for a subset of $2^{\omega}$
and every $p \in {\bold P}$,  if $p \forces \mu(\dot{A}) \leq \varepsilon
$ then
$$\mu^{\star}(\{x \in 2^{\omega} : \exists q \geq p \ q \forces
x \not \in \dot{A} \}) \geq 1-\varepsilon .$$

$\star_{3}[{\bold P}]$ iff
for every $A \subset {\bold V} \cap 2^{\omega}$ of positive measure
${\bold V}^{{\bold P}} \models \mu^{\star}(A) >0$. 

$\star_{4}[{\bold P}]$ 
iff for every sufficiently large cardinal $\kappa$,
and for every countable elementary submodel $N \prec H(\kappa,\in)$, if
${\bold P} \in N$ and $\{p_{n} : n \in \omega\} \in N$ is a
sequence of  ${\bold P}$ and
$\{\dot{A}_{n}: n \in \omega\} \in N$ is a sequence of elements of
${\bold P}$-names such that
for $n \in \omega$
$p_{n} \forces \dot{A}_{n}$ is a Borel set of measure $\leq \varepsilon_{n}$,
and $\lim_{n \rightarrow \infty} \varepsilon_{n}=0$ then for
every random real $x$ over $N$
there exists a condition $q \in {\bold P}$ such that
\begin{enumerate}
\item $q$ is $(N,{\bold P})$-generic,
\item $q \forces x$ is random over $N[G]$ for every ${\bold P}$-generic filter
over $N$ containing $p_{0}$,
\item there exists $n \in \omega$ such that $q \geq p_{n}$ and
$q \forces x \not \in \dot{A}_{n}$.
\end{enumerate}
\end{definition}

In \cite{JS} it is proved that
\begin{lemma}
For every notion of forcing ${\bold P}$,
\begin{enumerate}
\item If\/ ${\bold P}$ is weakly homogenous then $\star_{2}
[{\bold P}]
\leftrightarrow
\star_{3}
[{\bold P}]
$,
\item $\star_{1}
[{\bold P}]
\leftrightarrow
\star_{4}
[{\bold P}]
$.~$\QED$
\end{enumerate}
\end{lemma}

\begin{lemma}
Suppose that ${\bold P}$ has property $\star_{1}$. Then
${\bold V}^{{\bold P}}
\models$``${\bold V} \cap 2^{\omega}$ is not measurable''.
\end{lemma}
\Proof
It is enough to show that
${\bold V} \cap 2^{\omega}$ has positive outer measure.
Let $\{\dot{I}_{n} : n \in \omega\}$ be a ${\bold P}$-name for a
sequence of rational intervals such that
$p_{0} \forces \sum_{n \in \omega} \mu(\dot{I}_{n}) \leq \varepsilon < 1$.
Find sequences
$\{p_{n} : n \in \omega\}$,
$\{j_{n} : n \in \omega\}$,
and $\{I_{n} : n \in \omega\}$
such that for $n \in \omega$
\begin{enumerate}
\item $p_{n+1} \geq p_{n}$,
\item $p_{n+1} \forces \dot{I}_{j} = I_{j}$ for $j \leq j_{n}$,
\item $p_{n+1} \forces \sum_{j=j_{n}}^{\infty}
\mu(\dot{I}_{j}) \leq \varepsilon - \frac{1}{n}$.
\end{enumerate}
It is easy to see that
$\sum_{n \in \omega} \mu(I_{n}) \leq \varepsilon$.

Choose a countable, elementary submodel $N$ of $H(\kappa)$ containing
${\bold P}$ and
$\{p_{n},j_{n},\dot{I}_{n}, I_{n} : n \in \omega\}$.
Since $N$ is countable there exists $x \in {\bold V} \cap 2^{\omega}$
such that $x$ is a random real over $N$ and
$x \not \in \bigcup_{n \in \omega} I_{n}$.
Using $\star_{1}[{\bold P}]$ we get $q \geq p$ such that
$q \forces x \not \in \bigcup_{n \in \omega} \dot{I}_{n}$.

Since
$\{\dot{I}_{n} : n \in \omega\}$ was arbitrary it shows that
$${\bold V}^{{\bold P}}
\models \mu^{\star}({\bold V} \cap 2^{\omega}) =1 $$
which finishes the proof.~$\QED$

The lemma above would be even easier to prove if we assume
$\star_{3}[{\bold P}]$.
The reason for using property
$\star_{1}[{\bold P}]$ is in the following:
\begin{theorem}[\cite{JS}]
Suppose that $\{{\bold P}_{\xi}, \dot{{\bold Q}}_{\xi} :
\xi < \alpha\}$ is a countable support iteration such that
$\forces_{\xi}$ ``$\dot{{\bold Q}}_{\xi}$ has property
$\star_{1}$ for $\xi < \alpha$.
Let  ${\bold P}={\bold P}_{\alpha}$. Then ${{\bold P}}$  has property
$\star_{1}$. $\QED$
\end{theorem}
To construct the model satisfying $\ZFCa \ \& \ \neg w{\bold D}\ \&\
{\bold U}(c)\ \&\ \neg{\bold U}(m)\ \&\ \neg{\bold B}(m)$ we show that
forcing ${\bold Q}_{f,g}$ has property 
$\star_{1}$. At the first step we show that it has property
$\star_{3}$ i.e.
\begin{theorem}
Let $A \subset 2^{\omega}$ be such that $\mu(A) = \varepsilon_{0} >0$.
Then ${\bold V}^{{\bold Q}_{f,g}} \models \mu^{\star}(A)~>~0$.
\end{theorem}
\Proof
Suppose that this theorem is not true. Then there exists a set
$A \subset 2^{\omega}$  such that $\mu^{\star}(A) = \varepsilon_{0} >0$,
a condition $T \in {\bold Q}_{f,g}$
and a sequence $\{\dot{I}_{n} : n \in \omega\}$ of
${\bold Q}_{f,g}$-names for rational intervals such that
\begin{enumerate}
\item $T \forces \sum_{n=1}^{\infty} \mu(\dot{I}_{n}) =1$,
\item $T \forces A \subset \bigcap_{m \in \omega} \bigcup_{n > m} \dot{I}_{n}$.
\end{enumerate}
Let $s_{0}$ be the stem of $T$. By \ref{crucial} without losing generality
we can assume that there exists an increasing  sequence of natural numbers
$\{k_{n} : n \in \omega\}$ such that
\begin{enumerate}
\item For every $s \in T \cap \omega^{k_{n}} \ T^{[s]}$ forces a value
to $\{\dot{I}_{j} : j \leq n\}$,
\item $T \forces \sum_{n \geq \lh(s_{0})} \mu(\dot{I}_{n}) < \frac{1}{2} \cdot
\varepsilon_{0}$,
\item $\prod_{n=\lh(s_{0})}^{\infty} (1- \frac{1}{f(n)}) > \frac{1}{2}$.
\end{enumerate}
For $s \in T$ and $j \in \omega$ define
$$I^{s}_{j} = \left\{ \begin{array}{ll} I & \hbox{if }
T^{[s]} \forces \dot{I}_{j} = I\\
\emptyset & \hbox{otherwise}
\end{array} \right. .$$
Suppose that a function $h \in \omega^{\omega}$ witnesses that
$T \in {\bold Q}_{f,g}$ and consider a function $h' \in \omega^{\omega}$
such that $h'(n) \leq h(n)$ for $n \in \omega$.
\begin{claim}
For $x \in 2^{\omega}$ the following condition are equivalent:
\begin{enumerate}
\item There exists $T' \geq T$ such that $h'$ witnesses that
$T' \in {\bold Q}_{f,g}$ and $T' \forces x \not \in \bigcup_{n \in \omega}
\dot{I}_{n}$,
\item For every $k \geq \lh(s_{0})$ there exists a finite tree $t$ of height
$k$
such that
\begin{enumerate}
\item $t \subset T \cap \omega^{\leq k}$,
\item $|\suc_{t}(s)| \geq g(\lh(s),h'(\lh(s)))$ for $s \in t \cap \omega^{\geq
\lh(s_{0})}$,
\item If $s \in t \cap \omega^{k}$ then
$x \not \in \bigcup_{j \in \omega} I^{s}_{j}$.
\end{enumerate}
\end{enumerate}
\end{claim}
\Proof
$(1) \rightarrow (2)$
If $T'$ satisfies (1) then $T' \rest k$
satisfies (2)

$(2) \rightarrow (1)$ Build a sequence $\{t_{k} : k \in \omega\}$
satisfying (2) and apply the compactness theorem to construct $T'$.~$\QED$

Define a set $D \subset 2^{\omega}$ as follows:

$y \in D$ iff there exists $T' \in {\bold Q}_{f,g}$ such that
\begin{enumerate}
\item $T' \geq T$ has the same stem as $T$ (=$s_{0}$),
\item $T' \forces y \not \in \bigcup_{n \geq \lh(s_{0})} \dot{I}_{n}$,
\item $\forall n \geq \lh(s_{0}) \ \forall s \in T' \cap \omega^{n}
\ |\suc_{T'}(s)| \geq g(n,h(n)-1)$.
\end{enumerate}
Notice that the set $D$ is defined in ${\bold V}$ and since
$T \forces A \subset \bigcup_{n \geq \lh(s_{0})} \dot{I}_{n}$ we have
$\mu(2^{\omega} - D) > \varepsilon_{0}$.

For $ k \geq \lh(s_{0})$ define sets $D_{k}$ as follows:

$y \in D_{k}$ iff there exists a finite tree $t$ such that
\begin{enumerate}
\item $t \subset T \cap \omega^{\leq k}$ ,
\item $\forall n \geq \lh(s_{0}) \ \forall s \in t \cap \omega^{n}
\ |\suc_{t}(s)| \geq g(n,h(n)-1)$,
\item $\forall s \in t \cap \omega^{k} \ y \not \in \bigcup_{n \geq \lh(s_{0})}
I_{n}^{s}$.
\end{enumerate}
By the above claim $D = \bigcap_{k \in \omega} D_{k}$. Since sets
$D_{k}$ form a decreasing family we can find $k \in \omega$ such that
$\mu(2^{\omega} - D_{k}) > \varepsilon_{0}$.

For every $s \in T$ such that $\lh(s_{0}) \leq \lh(s) \leq k$ define
set $D_{k,s}$  as follows:

$y \in D_{k,s}$ iff there exists a finite tree $t$ such that
\begin{enumerate}
\item $t \subset T \cap \omega^{\leq k}$ and $t = t^{[s]}$,
\item $\forall n \geq \lh(s) \ \forall s' \in t \cap \omega^{n}
\ |\suc_{t}(s')| \geq g(n,h(n)-1)$,
\item $\forall s' \in t \cap \omega^{k} \ y \not \in \bigcup_{n \geq \lh(s_{0})}
I_{n}^{s'}$.
\end{enumerate}
Notice that $D_{k} = D_{k,s_{0}}$.
Observe also that for $s \in T \cap \omega^{k}$
$$\mu(2^{\omega} - D_{k,s}) \leq \sum_{n \geq \lh(s_{0})} \mu(I^{s}_{n}) <
\frac{\varepsilon_{0}}{2} .$$
\begin{claim}
Suppose that for some $m \in [\lh(s_{0}),k-1]$ and $s \in T \cap
\omega^{m}$, 
$$\mu(2^{\omega} - D_{k,t})~\leq~a \text{ for } t \in \suc_{T}(s).$$
Then
$$\mu(2^{\omega} - D_{k,s})
\leq \frac{a}{ 1-\displaystyle\frac{g(m,h(m)-1)}{g(m,h(m))}} .$$
\end{claim}
\Proof
Notice that
$y \not \in D_{k,s} $ iff $|\{t \in \suc_{T}(s) : y \not \in D_{k,t}\}| >
g(m,h(m)) - g(m,h(m)-1)$.
\begin{claim}
Let $N_{1}> N_{2}$ be two natural numbers. Suppose  that
$\{A_{j} : j \leq N_{1}\}$ is a family of subsets of $2^{\omega}$
of measure $\leq a$.
Let $U = \{ x \in 2^{\omega} : x$ belongs to at least $N_{2}$ sets
$A_{j} \}$. Then $$\mu(U) \leq a \cdot \frac{N_{1}}{N_{2}} .$$
\end{claim}
\Proof
Let $\chi_{A_{i}}$ be the characteristic function of the set $A_{i}$ for
$i \leq N_{1}$.
It follows that $\int \sum_{i \leq N_{1}} \chi_{A_{i}} \leq N_{1}
\cdot a$ and 
therefore
$$\mu\left(\left\{x \in 2^{\omega} :
\sum_{i \leq N_{1}} \chi_{A_{i}}(x)  \geq N_{2}\right\}\right) \leq
\frac{N_{1}}{N_{2}}\cdot 
a.~\QED$$

By applying the claim above we get
$$\mu(2^{\omega} - D_{k,s}) \leq a \cdot
\frac{g(m,h(m))}
{g(m,h(m))-g(m,h(m)-1)} =
\frac{a}{ 1-\displaystyle\frac{g(m,h(m)-1)}{g(m,h(m))}}  .~\QED$$

Finally by induction we have
$$\mu(2^{\omega} - D_{k}) =
\mu(2^{\omega} - D_{k,s_{0}}) \leq \frac{\varepsilon_{0}}{2} \cdot
\frac{1}{M}$$
where
$$M = \prod^{m=k}_{\lh(s_{0})}
\left(1-\frac{g(m,h(m)-1)}{g(m,h(m))}\right) \geq 
\prod^{m=k}_{\lh(s_{0})} \left(1-\frac{1}{f(m)}\right) > \frac{1}{2}.$$
Therefore
$\mu(2^{\omega} - D_{k}) < \varepsilon_{0}$ which gives a
contradiction.~$\QED$ 

Now we can prove
\begin{theorem}\label{star1}
${\bold Q}_{f,g}$ has property $\star_{1}$.
\end{theorem}
\Proof
We will need several definitions:

\begin{definition}
Let
$\{\dot{I}_{n}: n \in \omega\}$
be a ${\bold Q}_{f,g}$-name
for a sequence of rational intervals.
We say that $T \in {\bold Q}_{f,g}$ interprets
$\{\dot{I}_{n}: n \in \omega\}$
if there exists an increasing sequence
$\{k_{n}: n \in \omega\}$
such that
for every $j \leq n \in \omega$ and $s \in T \cap \omega^{k_{n}} \
T^{[s]}$ decides a value of $\dot{I}_{j}$ i.e.
$T^{[s]} \forces \dot{I}_{j} = I^{s}_{j}$ for some
rational interval $I^{s}_{j}$.
\end{definition}
By \ref{crucial} we know that
$$\{ T \in {\bold Q}_{f,g} : T \text{ interprets }
\{\dot{I}_{n}: n \in \omega\}\}$$  is  dense  in $ {\bold Q}_{f,g}$.
Suppose that $T \in {\bold Q}_{f,g}$. Subset
$S \subseteq T$ is called {\em front} if
for every branch $b$ through $T$ there exists $n \in \omega$ such
that $b \rest n \in S$.

Suppose that $D \subseteq {\bold Q}_{f,g}$ is an open  set.
Define
$$cl(D) = \{ T \in {\bold Q}_{f,g} :
\{s \in T: T^{[s]} \in D\} \ is \ a \ front \ in \ T\} .$$

Let
$\{\dot{I}_{n}: n \in \omega\}$
be a ${\bold Q}_{f,g}$-name for
a sequence of rational intervals
such that
for some $T_{0} \in {\bold Q}_{f,g} \ T_{0} \forces
\sum_{n=1}^{\infty} \mu(\dot{I}_{n}) < \varepsilon < 1$ and
$T_{0}$ interprets
$\{\dot{I}_{n}: n \in \omega\}$.

Let $N \prec H(\kappa)$ be a countable model containing
${\bold Q}_{f,g}, \ T_{0}, \ \{\dot{I}_{n}: n \in \omega\}$.

Define a set $Y \subseteq 2^{\omega}$ as follows:

$x \in Y$ iff there exists $\widehat{T} \in {\bold Q}_{f,g}$ such that
\begin{enumerate}
\item $\widehat{T} \leq T_{0}$,
\item If $D \in N$ is an open, dense subset of ${\bold Q}_{f,g}$
then there exists $T' \in cl(D) \cap N$ such that
$\widehat{T} \leq T'$,
\item $\widehat{T} \forces x \not \in \bigcup_{n \in \omega} \dot{I}_{n}$,
\item Suppose that
$J = \{\ddot{I}_{n} : n \in \omega\} \in N$
is a ${\bold Q}_{f,g}$-name for
a sequence of rational intervals
such that
$\forces
\sum_{n=1}^{\infty} \mu(\ddot{I}_{n}) < \infty$ and let
$D_{J} = \{
T \in {\bold Q}_{f,g}:  T$ interprets
$\{\ddot{I}_{n} : n \in \omega\} $
(with sequence $\{k_{n}^{T} : n \in \omega\}$).
Then there exists $T \in D_{J} \cap N$ and $k \in \omega$
such that

$\forall m \geq k \ \forall s \in \widehat{T} \cap \omega^{k_{m}^{T}} \ x \not \in
I^{T,s}_{m}$.
\end{enumerate}
Notice that (2) guarantees that $\widehat{T}$ is $(N,{\bold Q}_{f,g})$-generic
while (4) guarantees that $x$ is random over $N[G]$.

\begin{lemma}
\ \begin{enumerate}
\item $Y$ is a ${\bold \Sigma}^{1}_{1}$ set of reals (in ${\bold V}$),
\item $\mu(Y) \geq 1-\varepsilon$.
\end{enumerate}
\end{lemma}
\Proof
(1) It is easy to see that conditions (1)-(4) in the definition of
$Y$ are Borel provided that we  have an enumeration
(we can code as a real number) of the objects appearing in (2) and (4).

(2) easy computation using the fact that ${\bold Q}_{f,g}$ has property
$\star_{3}$ and $\star_{2}$.~$\QED$

Work in $N$. Let $G \subset \coll(\boldsymbol\aleph_{0}, 2^{\boldsymbol\aleph_{0}})$ be
generic over $N$ and let $x$ be a random real over $N[G]$.
Let ${\bold B}$ denotes the measure algebra.
Since parameters of the definition of $Y$ are in $N[G]$ 
we can
ask whether $N[G][x] \models x \in Y$.

Since in $N[G]$, $Y$ is a measurable set we can find two disjoint, Borel sets $A$
and
$B$ such that $\mu(A \cup B) =1$ and
$A \forces_{{\bold B}} x \in Y$ and $B \forces_{{\bold B}} x \not \in Y$.
Morover $\mu(A) \geq 1-\varepsilon$. In other words $A \subseteq Y \
a.e.$ and 
$B \subseteq 2^{\omega}-Y \ a.e.$

Since $x$ is a random real over $N$ as well we have
$$\coll(\boldsymbol\aleph_{0},2^{\boldsymbol\aleph_{0}}) \star {\bold B} \cong
Q_{x} \star \dot{{\bold R}} \cong
{\bold B} \star \dot{{\bold R}} $$
where $Q_{x}$ is the smallest subalgebra
which adds $x$.

Find a Borel set of positive measure $A^{\star}$ such that
$$N \models A^{\star} \forces_{{\bold B}}`` \exists p \in \dot{{\bold R}} \
p \forces x \in A''$$
and
$$N \models 2^{\omega}-A^{\star} \forces_{{\bold B}}
``\forces x \in B''.$$
It is clear that $A^{\star}-A$ has measure zero and therefore
$\mu(A^{\star}) \geq 1-\varepsilon$.

Notice that the  definitions above do not depend on the choice of
random real 
$x$ as long as $x \in A^{\star}$.
Thus if $x$ is {\em any} random real over $N$ such that $x \in A^{\star}$ then
we can find an $N$-generic filter $G \subset
\coll(\boldsymbol\aleph_{0},2^{\boldsymbol\aleph_{0}})$ 
such that $(G,x)$ is
$\coll(\boldsymbol\aleph_{0},2^{\boldsymbol\aleph_{0}}) \star {\bold B}$-generic  over $N$ and
$N[G][x] \models x \in Y$.
Since $Y$ is a ${\bold \Sigma}^{1}_{1}$ set it means that
${\bold V} \models x \in Y$.
In other words there exists a Borel set $A^{\star}$ of measure $\geq 1-
\varepsilon$ such that if $x \in {\bold V}\cap A^{\star}$ is a random real over
$N$
then $x \in Y$.

Now we finish the proof of the theorem. Let $N$ , $\{p_{n} : n \in
\omega\}$, 
$\{\dot{I}_{n} : n \in \omega\}$ and $x$  be such that
\begin{enumerate}
\item $p_{n+1} \geq p_{n}$ for $n \in \omega$,
\item $p_{n} \forces \dot{I}_{n} = I_{n}$ for $n \in \omega$,
\item $x \not \in \bigcup_{n \in \omega} I_{n}$,
\item $\sum_{n=1}^{\infty} \mu(I_{n}) = \varepsilon$.
\end{enumerate}
Define for $n \in \omega$, $Y_{n} =$ set $Y$ defined for
model $N$, condition $p_{n}$ and set $\{\dot{I}_{m+n} : m \in \omega\}$.

By the above remarks we can find Borel
sets $\{A^{\star}_{n} : n \in \omega\} \in N$
such that for $n \in \omega$
$\mu(A_{n}^{\star}) \geq 1-(\varepsilon-\sum_{j \leq n} \mu(I_{j}))$ and
for every $x \in {\bold V} \cap A_{n}^{\star}$ if
$x$ is random over $N$ then $x \in Y_{n}$.
Since  $\mu(\bigcup_{n \in \omega} A^{\star}_{n})=1$ if $x$ is random over
$N$ then $x \in A^{\star}_{n}$ for some $n \in \omega$.
Therefore $x \in Y_{n}$ and this finishes the proof as
$Y_{n} \subset Y_{0}$ for all $n \in \omega$.
From the fact that $x \in Y_{0}$ follows the existence of the condition
witnessing $\star_{1}$.~$\QED$

\begin{theorem}
$\Con(\ZFCa) \rightarrow
\Con(\ZFCa \ \& \ \neg w{\bold D}\ \&\
{\bold U}(c)\ \&\ \neg{\bold U}(m)\ \&\ \neg{\bold B}(m))$.
\end{theorem}
\Proof
Let  $\{{\bold P}_{\xi}, \dot{{\bold Q}}_{\xi} :
\xi < \boldsymbol\aleph_{2}\}$ be a countable support iteration such that
$\forces_{\xi}$ ``$\dot{{\bold Q}}_{\xi} = {\bold Q}_{f,g}$ for $\xi < \boldsymbol\aleph_{2}$.
Let ${\bold P} ={\bold P}_{\boldsymbol\aleph_{2}}$. Then
${\bold V}^{\bold P} \ \models \
\neg{\bold U}(m)$
because ${\bold P}$ has property $\star_{1}$ and
${\bold V}^{\bold P} \models
\neg{\bold B}(m)$ and $\neg w{\bold D}$
since ${\bold P}$ is $f$-bounding and $\omega^{\omega}$-bounding by
\ref{preserve} 
and \ref{bound}. Finally ${\bold V}^{\bold P} \models {\bold U}(c)$ by
\ref{2.1}.~$\QED$

\section{Rational perfect set forcing}
Our next goal is to construct a model for
$$\ZFCa \ \& \ w{\bold D}\
\&\ \neg{\bold D}\ \&\ {\bold U}(c) \ \&\ \neg{\bold U}(m)\ \&\
\neg{\bold 
  B}(m).$$ 
We will do it in the next section. H
This model is obtained as a $\omega_{2}$-iteration with countable
support of 
${\bold Q}_{f,g}$ and rational perfect set forcing. In this section we
will prove several facts about rational perfect set forcing which we
will need later.

Recall that rational perfect set forcing is defined
as follows:

$T \in {\bold R}$ iff $T$ is a perfect subtree of $\omega^{<\omega}$ and
for every $s \in T$ there exists $s \subseteq t \in T$ such that
$\suc_{T}(t)$ is infinite.

Elements of ${\bold R}$ are ordered by $\subseteq$.

Without loss of generality we can assume that for every $T \in {\bold
  R}$ and 
$s \in T$ the set $\suc_{T}(s)$ is either infinite or contains exactly
one element  since elements of this form are dense in ${\bold R}$.

For $T \in {\bold R}$ define
$$\spli(T) = \{s \in T : \suc_{T}(s) \hbox{ is   infinite } \} .$$
For $T,T' \in {\bold R}$ let

$T \geq_{0} T'$ if $T \geq T'$ and $T$ and $T'$ have the same stem.

$T' \geq_{n} T$ if $T' \geq T$ and for every $s \in \spli(T)$ if
exactly $n$ proper segments of $s$ belong to $\spli(T)$ then
$s \in \spli(T')$.

First we have to show that forcing ${\bold R}$ preserves outer measure.

\begin{definition}
Let
$\{\dot{I}_{n} : n \in \omega\}$
be an ${\bold R}$-name for
sequence of rational intervals such that
$\forces \sum_{n=1}^{\infty} \mu(\dot{I}_{n}) = \frac{1}{2}$.

We say that $T \in {\bold R}$ {\em interprets}
$\{\dot{I}_{n} : n \in \omega\}$
if for every $s \in \spli(T)$ there exist
rational intervals $\{I^{s}_{1}, \ldots , I^{s}_{n_{s}}\}$such that
\begin{enumerate}
\item $T^{[s]} \forces \forall j \leq n_{s} \ \dot{I}_{j} = I^{s}_{j}$,
\item for every $\varepsilon > 0$ and every branch $y$ through
$T$ there exists $m \in \omega$ such that for $k \geq m$
$$\mu\left(\bigcup_{j \leq n_{y \rest k}} I_{j}^{y \rest k}\right) \geq
\frac{1}{2} - \varepsilon .$$
\end{enumerate}
\end{definition}

\begin{lemma}\label{fusion}
Suppose that
$\{\dot{I}_{n} : n \in \omega\}$
is an  ${\bold R}$-name for
sequence of rational intervals. Assume that
$T \forces \sum_{n=1}^{\infty} \mu(\dot{I}_{n}) = \frac{1}{2}$.
Then there exists $\widehat{T} \geq T$ such that
$\widehat{T}$ interprets
$\{\dot{I}_{n} : n \in \omega\}$.
\end{lemma}
\Proof
Construct a sequence $\{T_{n} : n \in \omega\}\subset {\bold R}$ such that
$T_{n+1} \geq_{n} T_{n}$ for $n \in \omega$ as follows:

$T_{0} =T$ and suppose that $T_{n}$ is already constructed.

For every $s \in \spli(T_{n})$ such that exactly $n$ proper segments
of $s$ belong to $\spli(T_{n})$ and every $m \in \omega$ such that
$s^{\frown}\{m\} \in \suc_{T_{n}}(s)$ extend
$T^{[s^{\frown}\{m\}]}$ to decide a sufficiently long part of
$\{\dot{I}_{n} : n \in \omega\}$.
Paste all extensions together to get $T_{n+1}$.

Clearly $\widehat{T} = \bigcap_{n \in \omega} T_{n}$ has required
property.~$\QED$

Now we are ready to show:
\begin{theorem}
If $A \subseteq 2^{\omega}$ and $\mu(A)=1$ then
$\forces_{{\bold R}} \mu^{\star}(A) > 0$.
\end{theorem}
\Proof
Suppose not. Then there exists a measure one  set $A \subseteq
2^{\omega}$, 
a ${\bold R}$-name for
sequence of rational intervals
$\{\dot{I}_{n} : n \in \omega\}$
and a condition $T \in {\bold R}$ such that
\begin{enumerate}
\item $T \forces \sum_{n=1}^{\infty} \mu(\dot{I}_{n}) = \frac{1}{2}$.
\item $T \forces A \subset \bigcap_{n \in \omega} \bigcup_{m \geq n}
\dot{I}_{m}$.
\end{enumerate}
By the above lemma we can assume that
$T$ interprets
$\{\dot{I}_{n} : n \in \omega\}$.

For $s \in \spli(T)$ and $\varepsilon > 0$ define
$$h^{\varepsilon}(s) = \min\left\{j \in \omega: \sum_{i \leq j}
\mu(I^{s}_{i}) \geq 
\frac{1}{2}-\varepsilon\right\}$$
and
$$A^{\varepsilon}_{s} = \bigcup_{i \geq h^{\varepsilon}(s)} I^{s}_{i} .$$
Note that $h^\varepsilon(s)$ may be undefined for some $\varepsilon$
and $s$. 

Let $N$ be a countable, elementary submodel of $H(\kappa)$ for
sufficiently big $\kappa$.

Let $x \in A$ be a random real over $N$ .
The following holds in $N[x]$.
\begin{lemma}\label{4.3}
For every $\varepsilon > 0$ there exists a tree $T_{\varepsilon} \subset
T$ such that
\begin{enumerate}
\item $T_{\varepsilon}$ has no infinite branches.
\item for every $s \in T_{\varepsilon}$ either
$x \in A^{\varepsilon}_{s}$ or
$\{n \in \omega: s^{\frown}\{n\} \in \suc_{T}(s) - \suc_{T_{\varepsilon}}(s)\}$
is finite.
\end{enumerate}
\end{lemma}
\Proof
Fix $\varepsilon>0$.
For $s \in \spli(T)$ define an ordinal $r_{\varepsilon}(s)$ as follows:

$r_{\varepsilon}(s) =0 $ iff $x \in A^{\varepsilon}_{s}$,

$r_{\varepsilon}(s) = \limsup \{ r_{\varepsilon}(t)+1 : t \in \suc(s)
 $ and 
$r_{\varepsilon}(t)$ is defined$\}$.

In other words $r_{\varepsilon}(s)\geq \alpha$ iff for all $\beta <
\alpha$ 
there exists infinitely many $t \in \suc_T(s)$ such that
$r_{\varepsilon}(t) \geq \beta$. 

\begin{claim}
For every $s \in \spli(T)$ ordinal $r_{\varepsilon}(s)$ is well defined.
\end{claim}
\Proof
If not we inductively build a condition $T' \geq T^{[s]}$ such that
$r_{\varepsilon}(t)$ is not defined for all $t \in \spli(T')$. But then
$T' \forces x \not \in \bigcup_{n \geq h^{\varepsilon}(s)} \dot{I}_{n}$.
Contradiction.~$\QED$

Let $s_{0}$ be the stem of $T$. Define

$T_{\varepsilon} = \{ s \in T: s_{0} \subseteq s$ or
for all $k<l$ if
$r_{\varepsilon}(s\rest k)$ and
$r_{\varepsilon}(s\rest l)$
are defined then
$r_{\varepsilon}(s\rest k)>
r_{\varepsilon}(s\rest l))\}$

It is easy to see that  $T_{\varepsilon}$ has no  branches
since for every branch $y$ through $T$ there exists
$m \in \omega$ such that for $k \geq m \ r_{\varepsilon}(y \rest k)=0$.

On the other hand if $x \not \in A^{\varepsilon}_{s}$ then by the definition
of rank
the set
$\{t \in \suc(s) : r_{\varepsilon}(t) \geq r_{\varepsilon}(s)\}$ is at
most 
finite which verifies
(2).~$\QED$

By the above lemma for every $\varepsilon>0$ there exists
a tree $T_{\varepsilon}$ together with a function
$r_{\varepsilon} : \spli(T_{\varepsilon}) \longrightarrow \omega_{1}$
such that
$$\forall s,t \in \spli(T_{\varepsilon}) \ s \subset t \rightarrow
r_{\varepsilon}(s) > r_{\varepsilon}(t) .$$
Since $N[x]$ is a generic extension of $N$ there exists
Borel set $B \subset 2^{\omega}$  of positive measure such that
$$N \models B \forces_{{\bold B}} \forall \varepsilon > 0 \hbox{ there
  exist } 
r_{\varepsilon} \hbox{ and } T_{\varepsilon} \hbox{ as  in \ref{4.3} } .$$
Fix $\varepsilon_{0}= \mu(B)/{2}$ and let
$\dot{r}$ and $\dot{T}$ be
${\bold B}$-names for
$r_{\varepsilon_{0}}$ and  $T_{\varepsilon_{0}}$.

We can find Borel set $B' \subset B$ such that
$\mu(B') > \frac{1}{2} \cdot \mu(B)$ and for
$s \in \spli(T)$
\begin{enumerate}
\item
$\{n \in \omega:
\exists B'' \subset B' \ B'' \forces
s \in \dot{T} \ \& \ s^{\frown}\{n\} \not \in
\dot{T} \ \& \ \mu(B'' \cap A^{\varepsilon_{0}}_{s} )=0\}$ is finite,
\item $\{\alpha \in \omega_{1}:
\exists B'' \subset B' \ B'' \forces \dot{r}(s) = \alpha\}$ is finite.
\end{enumerate}
To show this we use the fact that the measure algebra ${\bold B}$ is
$\omega^{\omega}$-bounding and $\dot{T}$ is forced to satisfy
\ref{4.3}(2). 

Now define in $N$

$\widehat{T} = \{s \in T :
\exists B'' \subset B' \
B'' \forces_{{\bold B}} s \in \dot{T}\}$

and

$\widehat{r}(s) = \max(\{ \alpha < \omega_{1} :
\exists B'' \subset B' \
B'' \forces_{{\bold B}} s \in \dot{T}\ \& \ \dot{r}(s) =
\alpha\})$.

Notice that these definitions do not depend on the initial choice of
random real $x$ as long as $x \in B'$.
\begin{lemma}
\
\begin{enumerate}
\item $\widehat{T}$ is a subtree of $T$,
\item If $s \in \widehat{T}$ and $x \in B'$ is any random
real over $N$ such that $x \not \in A^{\varepsilon_{0}}_{s}$ and
$s \in \dot{T}[x]$ then
$\{n \in \omega: s^{\frown}\{n\} \in \dot{T}[x] -
\widehat{T}\}$ is finite,
\item If $t \subset s \in \widehat{T}$ then
$\widehat{r}(t)>
\widehat{r}(s)$.
\end{enumerate}
\end{lemma}
\Proof
(1) and (2) follow immediately from the definition of
$\widehat{T}$  and
the choice of the set $B'$.

(3) Suppose that $\widehat{r}(s) = \alpha$.
It means that there exists a set $B'' \subset B'$ such that
$$B'' \forces \dot{r}(s)=\alpha .$$
Thus
$$B'' \forces \dot{r}(t) \hbox{ is  well  defined  and } > \alpha
$$
so $\alpha < \widehat{r}(t)$.~$\QED$

In particular it follows from (3) that
the tree
$\widehat{T}$ is well-founded, i.e. has no infinite branches, and that
$\widehat{r} : \widehat{T} \longrightarrow \omega_{1}$
is a rank function such that
$$\forall s \subset t \in \widehat{T}\
\widehat{r}(s) >
\widehat{r}(t) .$$

By induction on rank  define sets $X_{s} \subset 2^{\omega}$
for $s \in \spli(\widehat{T})$ as follows:

If $\widehat{r}(s) = 0$ then
$X_{s} = A^{\varepsilon_{0}}_{s}$.
If $\widehat{r}(s) > 0$ then
$X_{s} = \{z \in 2^{\omega} : z$ belongs to all but finitely many
sets $X_{t}$ where $t$ is an immediate successor of $t$ is
$\spli(\widehat{T})\}$.

It is easy to check that
$\mu(X_{s}) \leq \varepsilon_{0}$ for $s \in \spli(\widehat{T})$.

Choose $x \in A \cap (B' - X_{s_{0}})$ which is random over $N$.
Since $x \not \in X_{s_{0}}$ we can find infinitely many
immediate successors $s$ of $s_{0}$ in $\spli(\widehat{T})$ such that
$x \not \in X_{s}$. Choose one of them, say $s_{1} \supset s_{0}$

such that $x \not \in X_{s_{1}}$ and $s_{1} \in \dot{T}[x]$.
By repeating this argument with
$s_{1}$ instead of $s_{0}$ and so on we construct a branch
through $\dot{T}[x]$. Contradiction  since the tree
$\dot{T}[x]$ is well-founded.~$\QED$.

By repeating the proof of \ref{star1} we get

\begin{theorem}
${\bold R}$ has property $\star_{1}$.~$\QED$
\end{theorem}

\section{Not adding dominating and Cohen reals}
In this section we construct models for
\begin{enumerate}
\item $\ZFCa \ \& \ {\bold D}\ \&\ {\bold B}(m)\ \&\ \neg{\bold B}(c)\ \&\ {\bold U}(m)$,
\item $\ZFCa \ \& \ w{\bold D}\ \&\ \neg{\bold D}\ \&\ \neg{\bold B}(c)\ \&\  {\bold
B}(m)\ \& \
{\bold U}(m)$,
\item $\ZFCa \ \& \ w{\bold D}\
\&\ \neg{\bold D}\ \&\ {\bold U}(c) \ \&\ \neg{\bold U}(m)\ \&\  \neg{\bold B}(m))$.
\end{enumerate}

We need the following definitions.

\begin{definition}
  Let ${\bold P}$ be a notion of forcing.
We say that ${\bold P}$ is {\em almost} $\omega^{\omega}$-{\em bounding} if
for every ${\bold P}$-name $\sigma$ such that $p \forces \sigma \in
\omega^{\omega}$
there exists a function $f \in {\bold V} \cap \omega^{\omega}$ such that
for every
subset $A \in {\bold V} \cap
[\omega]^{\omega}$ there exists $q \geq p$ such that
$$q \forces \exists^{\infty}n \in A \ \sigma(n) \leq f(n) .$$

We say that ${\bold P}$ is {\em weakly} $\omega^{\omega}$-{\em bounding} if
for every ${\bold P}$-name $\sigma$ such that $p \forces \sigma \in
\omega^{\omega}$
there exists a function $f \in {\bold V} \cap \omega^{\omega}$ such that
there exists $q \geq p$ such that
$$q \forces \exists^{\infty}n  \ \sigma(n) \leq f(n) .$$

\end{definition}

We will use the following two preservation theorems.
\begin{theorem}[\cite{Sh2}]\label{dominating}
Let $\{{\bold P}_{\xi}, \dot{{\bold Q}}_{\xi} : \xi < \alpha\}$ be a countable
support iteration such that for $\xi < \alpha$

$\forces_{\xi} \dot{{\bold Q}}_{\xi}$ is almost $\omega^{\omega}$-bounding.

Then ${\bold P}_{\alpha} = \lim_{\xi < \alpha} {\bold P}_{\xi}$ is
weakly $\omega^{\omega}$-bounding.~$\QED$ 
\end{theorem}

\begin{definition}
  Let ${\bold P}$ be a notion of forcing satisfying
axiom A. We say that ${\bold P}$ has {\em Laver property} if
there exists a function $f_{{\bold P}}
\in \omega^{\omega}$ such that
for every finite set $A \subset {\bold V}$, ${\bold P}$-name $\dot{a}$
, $p \in {\bold P}$ and $n \in \omega$ if $p \forces \dot{a} \in A$ then
there is $q \geq_{n} p$ and a set $B \subset A$ of size $\leq
f_{{\bold P}}(n)$
such that $q \forces \dot{a} \in B$.

\end{definition}

Notice that this definition is actually stronger than standard
definition of Laver property.

\begin{theorem}[\cite{JS1}]\label{cohen reals}
Let $S \subset \alpha$ and suppose that $\{{\bold P}_{\xi}, \dot{{\bold Q}}_{\xi} :
\xi < \alpha\}$ is a countable support iteration such that

$\forces_{\xi}$ ``$\dot{{\bold Q}}_{\xi}$ is a random real forcing'' if $\xi \in S$

$\forces_{\xi}$ ``$\dot{{\bold Q}}_{\xi}$ has Laver property'' if $\xi \not \in S$.

Let  ${\bold P}={\bold P}_{\alpha}$. Then no real in ${\bold V}^{{\bold P}}$ is
Cohen over ${\bold V}$. 
\end{theorem}
Now we can prove that:

\begin{theorem} \
\begin{enumerate}
\item $\Con(\ZFCa)\ \rightarrow \ \Con(\ZFCa\ \& \
{\bold D}\ \&\ {\bold B}(m)\ \&\ \neg{\bold B}(c)\ \&\ {\bold U}(m))$,
\item $\Con(\ZFCa)\ \rightarrow \ \Con(\ZFCa\ \& \
w{\bold D}\ \&\ \neg{\bold D}\ \&\ \neg{\bold B}(c)\ \&\  {\bold B}(m)\ \& \
{\bold U}(m))$,
\item $\Con(\ZFCa) \rightarrow
\Con(\ZFCa \ \& \ w{\bold D}\
\&\ \neg{\bold D}\ \&\ {\bold U}(c) \ \&\ \neg{\bold U}(m)\ \&\  \neg{\bold B}(m))$.
\end{enumerate}
\end{theorem}
\Proof
(1) Let $\{{\bold P}_{\xi}, \dot{{\bold Q}}_{\xi} :
\xi < \boldsymbol\aleph_{2}\}$ be a countable support iteration such that

$\forces_{\xi}$ ``$\dot{{\bold Q}}_{\xi}$ is a random real forcing'' if $\xi$ is even

$\forces_{\xi}$ ``$\dot{{\bold Q}}_{\xi}$ is Mathias forcing'' if $\xi$ is odd.

Let ${\bold P} ={\bold P}_{\boldsymbol\aleph_{2}}$. Then

${\bold V}^{\bold P} \ \models \
{\bold D}\ \&\ {\bold B}(m)\ \&\  {\bold U}(m)$
because Mathias and random reals are added cofinally in the iteration and

${\bold V}^{\bold P} \models
 \neg{\bold B}(c)$ by \ref{cohen reals}.

(2) Let $\{{\bold P}_{\xi}, \dot{{\bold Q}}_{\xi} :
\xi < \boldsymbol\aleph_{2}\}$ be a countable support iteration such that

$\forces_{\xi}$ ``$\dot{{\bold Q}}_{\xi}$ is a random real forcing'' if $\xi$ is even

$\forces_{\xi}$ ``$\dot{{\bold Q}}_{\xi}$ is Shelah forcing from \cite{BS}'' if $\xi$ is
odd.

Let ${\bold P} ={\bold P}_{\boldsymbol\aleph_{2}}$. Then

${\bold V}^{\bold P} \models \
w{\bold D}\  \&\          {\bold B}(m)\ \& \
{\bold U}(m)$ because
of properties of Shelah forcing and random forcing.
To  show that
${\bold V}^{\bold P} \models
\neg{\bold B}(c)$  we use  \ref{cohen reals} and the fact that Shelah
forcing has the Laver property.

(3) Let $\{{\bold P}_{\xi}, \dot{{\bold Q}}_{\xi} :
\xi < \boldsymbol\aleph_{2}\}$ be a countable support iteration such that

$\forces_{\xi}$ ``$\dot{{\bold Q}}_{\xi} \cong {\bold Q}_{f,g}$'' if $\xi$ is even

$\forces_{\xi}$ ``$\dot{{\bold Q}}_{\xi}\cong {\bold R}$'' if $\xi$ is odd.

Let ${\bold P} ={\bold P}_{\boldsymbol\aleph_{2}}$.
Since ${\bold R}$ is
has Laver property (\cite{M2}) exactly as in \ref{preserve} we show that
${\bold P}$ is $f$-bounding. Therefore
${\bold V}^{\bold P} \ \models \
\neg{\bold B}(m)$.
${\bold V}^{\bold P} \models
\neg {\bold U}(m)$ since
${\bold Q}_{f,g}$  and ${\bold R}$ have property $\star_{1}$. Also ${\bold
  V} \models w{\bold D} \ \& \ {\bold U}(c)$ since ${\bold R}$ adds
unbounded reals and by \ref{2.1}.

To finish the proof of (2) and (3) we have to check that
forcings used there do not add dominating reals.
By \ref{dominating} it is enough to verify that both Shelah forcing and
rational perfect set forcing are almost $\omega^{\omega}$-bounding and this will
be proved in the next theorem.~$\QED$

\begin{theorem}
\begin{enumerate}
\item Rational perfect set forcing ${\bold R}$ is almost
$\omega^{\omega}$-bounding,
\item The Shelah forcing is almost $\omega^{\omega}$-bounding.
\end{enumerate}
\end{theorem}
\Proof
Let $\sigma$ be an ${\bold R}$-name such that $T \forces \sigma \in
\omega^{\omega}$ for some $T \in {\bold R}$.
As in \ref{fusion} we can assume that
for every $s \in \spli(T)$ and $ t \in \suc_T(s)$, $T^{[t]}$ decides
the  value of $\sigma \rest \lh(s)$.
Notice that in this case every branch through $T$ gives an interpretation
to $\sigma$.
Let $N$ be a countable, elementary submodel of $H(\kappa)$ such that
${\bold R},\ T$ and $\sigma$ belong to $N$.
Let $g \in {\bold V} \cap \omega^{\omega}$ be a function which dominates
all elements of $N \cap \omega^{\omega}$.
Fix a set $A \in {\bold V} \cap [\omega]^{\omega}$.
Since forcing ${\bold R}$ has absolute definition it is enough to show that
for every $m \in \omega$ and every condition $ T' \in N \cap
{\bold R}$, $T \leq T'$
there exists a condition $T'' \in N \cap {\bold R}$, $T' \leq T''$  and
$n \in A - [0,m]$ such that
$N \models T'' \forces \sigma(n) \leq g(n)$.
Choose $T' \geq T$ and let $b \in N$ be a branch through $T'$.
Let $\sigma_{b} \in N \cap \omega^{\omega}$ be the interpretation of
$\sigma$ obtained using $b$. By the assumption there exists
$n \in A,\ n \geq m$  such that $\sigma_{b}(n) \leq g(n)$. Choose
$T'' = T^{'[t]}$ where $t = b \rest n$.

(2) The proof presented here uses notation from \cite{BS}. Since the definition
of Shelah's forcing and all the necessary lemmas can be found in \cite{BS}
we give here only a skeleton of the proof.

Let $p=(w,T) \in {\bold S}$ and let $\tau$ be an ${\bold S}$-name for
an element of $\omega^{\omega}$. Let $q$ be a pure extension
of $p$ satisfying 2.4 of \cite{BS}.
Suppose that $q=(w,t_{0},t_{1}, \ldots)$.
We define by induction a sequence $\{q_{l} : l \in \omega\}$
satisfying the following conditions:
\begin{enumerate}
\item $q_{0}=q$,
\item $q_{l+1} = (w, t_{0}^{l+1}, t_{1}^{l+1} , \ldots)$
is  an $l$-extension of $q_{l}$,
\item if $k \leq l+1$ and $(w,w') \in t_{0}^{l+1} \ldots t_{k}^{l+1}$
and $w' \cap [n(t^{l+1}_{k}), m(t^{l+1}_{k})) \neq \emptyset$ when
$t^{l+1}_{k} \in K_{n(t^{l+1}_{k}),m(t^{l+1}_{k})}$ then
$(w', t^{l+1}_{k+1}, t^{l+1}_{k+1}, \ldots)$ forces value for
$\tau \rest k$,
\item $Dp(t^{l+1}_{l+1}) > l$.
\end{enumerate}
Before we construct this sequence let us see that this is enough to
finish the proof.

Let $q^{\star} = (w, t^{1}_{1},t^{2}_{2}, \ldots )$.
By $(4)$, $q^{\star} \in {\bold S}$.

Let $g(n) = \max \{k : \exists w' \ (w,w') \in t^{1}_{1} \ldots t^{n}_{n}$
and
$(w', t^{n+1}_{n+1}, t^{n+2}_{n+2}, \ldots ) \forces \tau(n) = k\}$
for $n \in \omega$.

Clearly $g \in \omega^{\omega}$.
Suppose that $A \subset \omega$. Define
$$p_{A} = (w, (t^{i}_{i} : i \in A)) .$$
It is easy to see that
$$p_{A} \forces \exists^{\infty}n \in A \ \tau(n) \leq g(n) $$
which finishes the proof.

We build the sequence $\{q_{l}: l \in \omega\}$ by induction on $l$.
Suppose that $q_{l}$ is already given.
By the definition of ${\bold S}$ it is enough to build the condition for
some fixed $w^{\star}= w \cap m(t^{l}_{0}, \ldots t^{l}_{l})$.

Define a function $C: \omega^{<\omega} \longrightarrow 2$ as follows:

$$C(v) = 1 \hbox{ iff } \exists k \ (w^{\star},v) \in t^{l}_{l+1}, \ldots,
t^{l}_{k} \hbox{ and } (v, t^{l}_{k+1}, t^{l}_{k+2} , \ldots )
\hbox{ forces value for } \tau(l) .$$
Using lemma 2.6 from \cite{BS} we get a condition where the function 
$C$ is constantly $0$ or $1$. The first is impossible since the set of
conditions forcing a value for $\tau(l)$ is dense. 
Therefore we get a condition $q=q_{l+1}$ on which $C$ is constantly $1$. 
Moreover we can assume that $q_{l+1}$ is an $l$-extension of $q_{l}$. 

This finishes the induction and the proof.~$\QED$

\end{document}